# STOCHASTIC CALCULUS FOR FRACTIONAL BROWNIAN MOTION WITH HURST EXPONENT $H > \frac{1}{4}$: A ROUGH PATH METHOD BY ANALYTIC EXTENSION

BY JÉRÉMIE UNTERBERGER

*Université Henri Poincaré Nancy I*

The $d$-dimensional fractional Brownian motion (FBM for short) $B_t = ((B_t^{(1)}, \ldots, B_t^{(d)}), t \in \mathbb{R})$ with Hurst exponent $\alpha$, $\alpha \in (0, 1)$, is a $d$-dimensional centered, self-similar Gaussian process with covariance $\mathbb{E}[B_s^{(i)} B_t^{(j)}] = \frac{1}{2}\delta_{i,j}(|s|^{2\alpha} + |t|^{2\alpha} - |t-s|^{2\alpha})$. The long-standing problem of defining a stochastic integration with respect to FBM (and the related problem of solving stochastic differential equations driven by FBM) has been addressed successfully by several different methods, although in each case with a restriction on the range of either $d$ or $\alpha$. The case $\alpha = \frac{1}{2}$ corresponds to the usual stochastic integration with respect to Brownian motion, while most computations become singular when $\alpha$ gets under various threshhold values, due to the growing irregularity of the trajectories as $\alpha \to 0$.

We provide here a new method valid for any $d$ and for $\alpha > \frac{1}{4}$ by constructing an approximation $\Gamma(\varepsilon)_t$, $\varepsilon \to 0$, of FBM which allows to define iterated integrals, and then applying the geometric rough path theory. The approximation relies on the definition of an analytic process $\Gamma_z$ on the cut plane $z \in \mathbb{C} \setminus \mathbb{R}$ of which FBM appears to be a boundary value, and allows to understand very precisely the well-known (see [5]) but as yet a little mysterious divergence of Lévy's area for $\alpha \to \frac{1}{4}$.

**0. Introduction.** The (two-sided) fractional Brownian motion $t \to B_t$, $t \in \mathbb{R}$ (FBM for short) with Hurst exponent $\alpha$, $\alpha \in (0, 1)$, defined as the centered Gaussian process with covariance

$$(0.1) \qquad \mathbb{E}[B_s B_t] = \tfrac{1}{2}(|s|^{2\alpha} + |t|^{2\alpha} - |t-s|^{2\alpha}),$$

is a natural generalization in the class of Gaussian processes of the usual Brownian motion, in the sense that it exhibits two fundamental properties shared with Brownian motion, namely, it has stationary increments, namely $\mathbb{E}[(B_t - B_s)(B_u - B_v)] = \mathbb{E}[(B_{t+a} - B_{s+a})(B_{u+a} - B_{v+a})]$ for every $a, s, t, u, v \in \mathbb{R}$, and it is self-similar, namely

$$(0.2) \qquad \forall \lambda > 0 \qquad (B_{\lambda t}, t \in \mathbb{R}) \stackrel{\text{(law)}}{=} (\lambda^\alpha B_t, t \in \mathbb{R}).$$

One may also define a $d$-dimensional vector Gaussian process (called: $d$-*dimensional fractional Brownian motion*) by setting $\mathbf{B}_t = (B_t^{(1)}, \ldots, B_t^{(d)})$, where $(B_t^{(i)}, t \in \mathbb{R})_{i=1,\ldots,d}$ are $d$ independent (scalar) fractional Brownian motions.

Its theoretical interest lies, in particular, in the fact that it is (up to normalization) the only Gaussian process satisfying these two properties.

A standard application of Kolmogorov's theorem shows that FBM has a version with $(\alpha - \varepsilon)$-Hölder paths for every $\varepsilon > 0$. In particular, all its paths possess finite $q$-variation for every $q > \frac{1}{\alpha}$, in the sense that

$$(0.3) \qquad \sup_{s = t_0 < \cdots < t_n = t} \left( \sum_{l=0}^n |B_{t_l} - B_{t_{l-1}}|^q \right) < \infty \qquad \text{a.s.,}$$

where the sum ranges over all partitions $(s = t_0 < t_1 < \cdots < t_n = t)$ of any order $n$ of the interval $[s, t]$.

There has been a widespread interest during the past ten years in constructing a stochastic integration theory with respect to FBM and solving stochastic differential equations driven by FBM. One of the strategies consists in trying to imitate the Brownian (or, more generally, martingale or even semi-martingale) construction by defining, for instance, the so-called *symmetric integral*

$$\int_0^t X_s \, dB_s = w - \lim_{\varepsilon \to 0} \int_0^t \frac{X_s + X_{s+\varepsilon}}{2} \frac{B_{s+\varepsilon} - B_s}{2\varepsilon} \, ds.$$

This idea and generalizations thereof (by using the classical Newton–Côtes approximation scheme for integration) led to the so-called Russo–Vallois theory of integration (see, e.g., [9, 18, 19]) which is well defined for the *scalar* fractional Brownian motion with any Hurst exponent $\alpha > 0$.

Another approach by using the Malliavin calculus and Skorokhod integration is due to Decreusefond and Üstünel (see [6]) and Cheridito and Nualart (see [4]). Again in the scalar case, they construct an extension of the divergence operator for every Hurst exponent $\alpha > 0$ and retrieve, in particular, that $\alpha = \frac{1}{6}$ is a barrier for the existence of the symmetric integral.

A very different approach, more suitable for nonscalar FBM, is ongoing since the work of Lyons; see [14]. This is the one we are going to follow. We shall need to recall the main results and notation of [5, 14] before we can state our own results.

Let $V$ be a Euclidean space, and $V^{\otimes k} = V \otimes \cdots \otimes V$ be the $k$th tensor product of $V$ endowed with a compatible norm $\|\cdot\|$. Lyons starts by defining a space of nonsmooth paths $\boldsymbol{\Gamma}:\mathbb{R}_+ \to V$ on V, called *geometric rough paths*, by completing the space $\Omega^\infty(V)$ of smooth paths. For $\Gamma \in \Omega^\infty(V)$, $k \geq 1$, and $s, t \geq 0$, we let

$$(0.4) \qquad \boldsymbol{\Gamma}^k_{s,t} = \int_s^t \boldsymbol{\Gamma}'_{t_1} \, dt_1 \otimes \int_s^{t_1} \boldsymbol{\Gamma}'_{t_2} \, dt_2 \otimes \cdots \otimes \int_s^{t_{k-1}} \boldsymbol{\Gamma}'_{t_k} \, dt_k.$$

This is an element of $V^{\otimes k}$, with components $(\Gamma^k_{s,t})_{i_1,\ldots,i_k}$ defined in an obvious way by projecting $\boldsymbol{\Gamma}'_{t_1}, \ldots, \boldsymbol{\Gamma}'_{t_k}$ onto the corresponding components. Repeating the same construction for another smooth path $\boldsymbol{\Delta}$, one may define the following *q-variation* distance for every $q > 1$:

$$(0.5) \qquad \delta_q(\boldsymbol{\Gamma}, \boldsymbol{\Delta})_{s,t} = \sup_{u \in [s,t]} \|\boldsymbol{\Delta}_u - \boldsymbol{\Gamma}_u\| + \sup_{1 \leq k \leq [q]} d_k(\boldsymbol{\Gamma}, \boldsymbol{\Delta})_{s,t},$$

where $[q]$ is the entire part of $q$ and $d_k$ is the distance defined by

$$(0.6) \qquad d_k(\boldsymbol{\Gamma}, \boldsymbol{\Delta})_{s,t} = \left( \sup_{s=t_0 < \cdots < t_n = t} \sum_{l=1}^n \|\boldsymbol{\Delta}^k_{t_{l-1},t_l} - \boldsymbol{\Gamma}^k_{t_{l-1},t_l}\|^{q/k} \right)^{k/q}.$$

Note that $d_k$ is the distance associated to a semi-norm that we shall call the $d_k$-*norm* (which depends on the interval $[s,t]$). The completion of $\Omega^\infty(V)$ with respect to the $q$-variation distance will be denoted by $\Omega_q(V)$.

The following extension principle is proved in [14]; see also the book by Lyons and Qian [15] or the nice review [12].

EXTENSION PRINCIPLE. For a given smooth function $f : \mathbb{R} \to \mathbb{R}^n \to \mathrm{Lin}(\mathbb{R}^d, \mathbb{R}^n)$, let $Y$ be the solution of the stochastic differential equation $dY_t = f(t, Y_t) \, dX_t$. Then the Itô map

$$(0.7) \qquad \mathcal{I} : \Omega^\infty(\mathbb{R}^d) \to \Omega^\infty(\mathbb{R}^n), \qquad X \to Y$$

is continuous with respect to the $q$-variation distance. Hence, $\mathcal{I}$ admits a unique continuous extension $\widetilde{\mathcal{I}} : \Omega^q(\mathbb{R}^d) \to \Omega^q(\mathbb{R}^n)$.

This result (obtained by the usual Picard iteration method) contains, in particular, an extension theorem for iterated integrals $X \to \int_0^t dX_{s_1} f_1(X_{s_1}) \times \int_0^{s_1} \cdots \int_0^{s_{n-1}} dX_{s_n} f_n(X_{s_n})$.

Hence, the problem of solving differential equations driven by FBM or giving a sense to iterated integrals with respect to FBM may be reduced to that of constructing an approximation of FBM which converges in the sense of rough paths.

The paper is organized as follows. We give in the first section a series decomposition of the fractional Brownian motion (see Lemma 1.1, Definition 1.2 and Corollary 1.7),

$$(0.8) \qquad \Gamma_t = \sum_{k \geq 0} \left( \int_0^t f_k(s) \, ds \right) \xi_k^+ + \sum_{k \geq 0} \left( \overline{\int_0^t f_k(s) \, ds} \right) \xi_k^-,$$

where $\xi_k^+, \xi_k^- = \overline{\xi_k^+}$ are independent complex normal Gaussian random variables, and show that $\Gamma_t$ may be seen as the restriction to the real axis of the real part of a Gaussian process $\Gamma_z^+$ defined on the closed upper half-plane $\bar{\Pi}^+$ which is a.s. analytic in $z$ on $\Pi^+$. Equivalently (think of Schwarz's reflection lemma), one may extend $\Gamma_z^+$ to the lower half-plane by simply taking its conjugate, $\Gamma_{\bar{z}}^- := \overline{(\Gamma_z^+)}$. In this way, FBM appears as the boundary value of a Gaussian process defined on $\mathbb{C} \setminus \mathbb{R}$. Lemma 1.5 is important for it shows that the $(\alpha - \varepsilon)$-Hölder property of FBM extends to the process on the closed upper-half plane.

Interestingly enough, the $\Gamma$-process may be seen (after a Cayley transform) as a *random entire series* living in the unit disk. Lemma 1.5 implies then that convergence to FBM on the unit circle is uniform in all directions (even tangentially).

The reader may be puzzled by the fact that we actually obtain a decomposition of the *fractional noise*, that is, the (distribution-valued) derivative of FBM (hence the integrals over $s$ in the above definition of $\Gamma_t$). This is directly related to the *stationarity* (and simplicity of course!) of the "covariance kernel" which is (up to a coefficient) the function $(x,y) \to |x - y|^{2\alpha - 2}$ (at least for $\alpha > \frac{1}{2}$, for otherwise it is singular). This kernel has many known decompositions in terms of orthogonal polynomials or special functions; the choice we made for the functions $f_k$ is the simplest one, but there are other natural possibilities which also lead to a straightforward analytic extension to the upper half-plane. This approach may be explained in "physical" terms by saying that we approximate the singular kernel $(x - y)^{2\alpha - 2}$ by a regular kernel $((x - y) \pm i\varepsilon)^{2\alpha - 2}$—or, more precisely, $(\pm i(x - y) + \varepsilon)^{2\alpha - 2}$ (see below)—where $\varepsilon$ is a short-distance cut-off.

The paper is organized as follows.

Section 1 is devoted to the construction of the $\Gamma$-process, defined as the above series (0.8), and to the proof that its boundary value on $\mathbb{R}$ has the same law as FBM. The main result of the section is Lemma 1.5.

Section 2 can be skipped for a first reading since the results will not be used in the following sections. It is devoted to the proof of convergence theorems for the series. It may be of interest for those interested in simulating FBM. Unfortunately, it is absolutely of no use when trying to construct a stochastic integration with respect to FBM, except maybe in the easy case $\alpha > \frac{1}{2}$ where simpler methods are more adapted.

In Section 3 we introduce our $\varepsilon$-approximation $\Gamma(\varepsilon)_t$ of FBM. The process is simply defined as the real part of $\Gamma_{t+i\varepsilon}^+$. The $\varepsilon$-shift of the variable $t$ into the upper half-plane transforms FBM into an analytic process. Contrary to the usual schemes of linear interpolation, *exact* computations are still possible when replacing FBM with $\Gamma(\varepsilon)_t$.

Section 4 is the hard part. It is devoted to the proof of the existence of the rough path limit of $\Gamma(\varepsilon)$ for $\varepsilon \to 0$ when $\alpha > \frac{1}{4}$, which implies the possibility to define stochastic integration with respect to FBM and to solve FBM-driven stochastic differential equations. It is mostly interesting from a practical point of view, since the existence of a rough path limit for $\alpha > \frac{1}{4}$ is already known from the paper by Coutin and Qian [5]. Namely, we compute explicitly the second moment of the Lévy area of $\Gamma(\varepsilon)$ in terms of hypergeometric functions (see Theorem 4.4), from which the divergent terms in the limit $\varepsilon \to 0$ (for $\alpha < \frac{1}{4}$) may be precisely identified. We also compute the second moment of the Lévy volume of $\Gamma(\varepsilon)$ for a three-dimensional FBM. It is less explicit, given as a sum of terms, each term being an integrated hypergeometric function with respect to a power kernel (see Section 4.3) which converges either for $\alpha > \frac{1}{4}$ or for $\alpha > \frac{1}{6}$.

It would be interesting to know whether the piecewise linear approximation used in [5] and the analytic approximation defined in this paper define in the limit the same stochastic integration theory and give the same solution to stochastic differential equations. We shall come back to this at the end of Section 4.1.

Note finally that Feyel and de la Pradelle (see [8]) proved an Itô formula for FBM (actually, for the Liouville fractional Brownian motion, defined as fractional integrated Brownian motion) by using an analytic continuation *in the parameter* $\alpha$ for $\operatorname{Re} \alpha > \frac{1}{2}$. Their approach is totally different, so are the results.

*Notation.* We will use the following notation throughout the article.

The Hurst exponent $\alpha \in (0,1)$ is fixed once and for all. It is intended implicitly to differ from $\frac{1}{2}, \frac{1}{4}, \frac{1}{6}, \ldots$ unless explicitly mentioned. One could actually suppose that $\alpha < \frac{1}{2}$ since there are no singularities at all for $\alpha > \frac{1}{2}$ (and, hence, the procedure of analytic extension to the upper half-plane is pointless), and the theory of stochastic integration of Brownian motion (corresponding to the case $\alpha = \frac{1}{2}$) is pretty well understood without all this machinery. Nevertheless, all computations below also hold true for $\alpha > \frac{1}{2}$ (unless explicitly stated).

The one-dimensional fractional Brownian motion process with Hurst exponent $\alpha$ will be denoted by $(B_t, t \in \mathbb{R})$. A $d$-dimensional FBM with independent exponents will be written as $\mathbf{B} = (B^{(1)}, \ldots, B^{(d)})$. Alternatively, a process with the same law as $(B_t, t \in \mathbb{R})$ may be written as $(\Gamma_t, t \in \mathbb{R})$ if necessary. The process $\Gamma$ will be obtained as the boundary value of an analytic process defined on $\mathbb{C} \setminus \mathbb{R}$.

If $x \in \mathbb{R}$ and $k = 0, 1 \ldots$, then $(x)_k = x(x+1) \cdots (x+k-1) = \frac{\Gamma(x+k)}{\Gamma(x)}$ is the Pochhammer symbol.

If $x \in \mathbb{R}$, then $E(x) = \sup\{n \in \mathbb{Z} | n \le x\}$ is the entire part of $x$.

If $z \in \mathbb{C} \setminus \mathbb{R}_-$, and $\beta \in \mathbb{C}$, then $z^\beta := \exp \beta \ln(z)$ is defined by using the usual determination of the logarithm.

If $x \in \mathbb{R}$, $x \neq 0$, then $\mathrm{sgn}(x) = \frac{x}{|x|}$ is the sign of $x$.

The Poincaré upper half-plane is denoted by $\Pi^+ = \{x + iy \in \mathbb{C} | y > 0\}$, its closure by $\bar{\Pi}^+ = \{x + iy \in \mathbb{C} | y \geq 0\}$. Similarly, $\Pi^-$ and $\bar{\Pi}^-$ denote the open, respectively, closed lower half-planes.

In the next proposition, we summarize for the unfamiliar reader the few things that should be known on the Cayley transform.

PROPOSITION 0.1.  *The Cayley transform* $\mathcal{C}: t \to z = \frac{t-i}{t+i}$ *is a bijection:*

(i) *from the one-point compactification $\overline{\mathbb{R}}$ of the real line onto the unit circle;*

(ii) *from the set $\{t \in \mathbb{C} | \mathrm{Im}\, t \geq 0\} \amalg \{\infty\}$ (equal to the one-point compactification of $\bar{\Pi}^+$) onto the closed unit disk $\bar{D}$.*

*In particular,* $\mathcal{C}(0) = -1$ *and* $\mathcal{C}(\infty) = 1$. *The inverse transform reads* $t = \mathcal{C}^{-1}(z) = i\frac{1+z}{1-z}$.

**1. Definition of the $\Gamma$-process and connection with the first chaos of FBM.**
For $\alpha > \frac{1}{2}$, the first chaos of FBM may be constructed by defining $\mathcal{H}_\alpha$ to be the Hilbert space completion of the space of step functions with respect to the scalar product

$$(1.1) \qquad \langle f, g \rangle_{\mathcal{H}_\alpha} = \alpha(2\alpha - 1) \int_{-\infty}^{\infty} \int_{-\infty}^{\infty} f(x)\overline{g(y)} |x - y|^{2\alpha - 2} \, dx \, dy$$

and setting

$$\mathbb{E}\left[\left(\int f(t)\, dB_t\right)\left(\int g(s)\, dB_s\right)\right] = \langle f, g \rangle_{\mathcal{H}_\alpha}$$

for $f, g \in \mathcal{H}_\alpha$. The space $\mathcal{H}_\alpha$ may be called a *reproducing kernel Hilbert space* for FBM. Note that the convolution kernel $K(x, y) = |x - y|^{2\alpha - 2}$ is in $L^1_{loc}$ only in the range $\alpha > \frac{1}{2}$. If $\alpha < \frac{1}{2}$, then this definition does not make sense any more and the construction of a reproducing kernel Hilbert space is more complicated, relying, for instance, on a mapping to $L^2(\mathbb{R})$ (and hence to Brownian motion) by means of fractional derivatives; see [17].

One would like to define the "derivative" of FBM as "the" process with covariance $\mathbb{E}[B'_x B'_y] = \alpha(2\alpha - 1)|x - y|^{2\alpha - 2}$, but this has no meaning at all because of the nonintegrability of the singularity on the diagonal $x = y$. Note also the awkwardness of the absolute value in the kernel $|x - y|^{2\alpha - 2}$, which is, however, necessary since the power functions are multivalued.

By going over to the complex plane, all these problems may be avoided, as we shall show presently.

LEMMA 1.1. *Let $f_k, k \geq 0$ be the analytic functions defined on the upper half-plane as*

$$f_k(z) = 2^{\alpha-1}\sqrt{\frac{\alpha(1-2\alpha)}{2\cos\pi\alpha}}\sqrt{\frac{(2-2\alpha)_k}{k!}}\left(\frac{z+i}{2i}\right)^{2\alpha-2}\left(\frac{z-i}{z+i}\right)^k.$$

*If $z, w \in \Pi^+$, then the series $\sum_{k\geq 0} f_k(z)\overline{f_k(w)}$ converges in absolute value and*

$$\sum_{k\geq 0} f_k(z)\overline{f_k(w)} = \frac{\alpha(1-2\alpha)}{2\cos\pi\alpha}(-i(z-\bar{w}))^{2\alpha-2}.$$

REMARK. The above definition makes sense since $\text{Re}(-i(z-\bar{w})) > 0$ (recall from the Introduction that the fractional powers are defined by means of the usual determination of the logarithm on $\mathbb{C} \setminus \mathbb{R}_-$).

PROOF OF LEMMA 1.1. By Proposition 0.1 in the Introduction, the Cayley transform $z \to \frac{z-i}{z+i}$ maps the upper half-plane onto the unit disk, hence, the convergence of the series $S(z,\bar{w}) := \frac{2\cos\pi\alpha}{\alpha(1-2\alpha)}\sum_{k\geq 0} f_k(z)\overline{f_k(w)}$. An explicit computation yields

$$S(z,\bar{w}) = 2^{2\alpha-2}\left(\frac{z+i}{2i}\right)^{2\alpha-2}\left(\frac{\bar{w}-i}{-2i}\right)^{2\alpha-2}\sum_{k\geq 0}\frac{(2-2\alpha)_k}{k!}\left(\frac{z-i}{z+i}\frac{\bar{w}+i}{\bar{w}-i}\right)^k$$

(1.2) $$= 2^{2\alpha-2}\left(\frac{z+i}{2i}\right)^{2\alpha-2}\left(\frac{\bar{w}-i}{-2i}\right)^{2\alpha-2}\left(1-\frac{z-i}{z+i}\frac{\bar{w}+i}{\bar{w}-i}\right)^{2\alpha-2}$$

$$= (-i(z-\bar{w}))^{2\alpha-2}. \qquad \square$$

DEFINITION 1.2. Let $\xi_k^1, \xi_k^2$, $k = 0, 1, \ldots$ be independent centered real Gaussian random variables with variance 1, and let $\xi_k^\pm = \xi_k^1 \pm i\xi_k^2$. The series

$$\sum_{k\geq 0} f_k(z)\xi_k^+ \qquad (z \in \Pi^+)$$

defines a random process $\Gamma_z'^+$; we denote by $\Gamma_{\bar{z}}'^-$ its complex conjugate $\sum_{k\geq 0} \overline{f_k(z)}\xi_k^-$.

With this notation, $\Gamma'^+$ (resp. $\Gamma'^-$) lives on $\Pi^+$ (resp. $\Pi^-$). The above lemma yields immediately the covariance

$$\mathbb{E}[\Gamma'^+(z)\Gamma'^+(w)] = \mathbb{E}[\Gamma'^-(\bar{z})\Gamma'^-(\bar{w})] = 0,$$

$$\mathbb{E}[\Gamma'^+(z)\Gamma'^-(\bar{w})] = \frac{\alpha(1-2\alpha)}{2\cos\pi\alpha}(-i(z-\bar{w}))^{2\alpha-2}.$$

By using the Cayley transform $z \to \zeta = \frac{z-i}{z+i}$, the process $\Gamma'^+$ may be rewritten (up to a prefactor) as a random entire series[1], namely,

$$(1.3) \quad \Gamma'^+_{z(\zeta)} = \left(\frac{1}{1-\zeta}\right)^{2\alpha-2} 2^{\alpha-1} \sqrt{\frac{\alpha(1-2\alpha)}{2\cos\pi\alpha}} \sum_{k\geq 0} \sqrt{\frac{(2-2\alpha)_k}{k!}} \zeta^k \cdot \xi_k^+.$$

Random series have been extensively studied; see, for instance, the book [10] by Kahane. The facts we need on random entire series are elementary, so we chose to give self-contained proofs. The next lemma will allow us to prove that $\Gamma'^+$ is analytic in the upper half-plane.

LEMMA 1.3.

1. Let $(a_n)_{n=0,1,\ldots}$ be independent random variables. Then there exists $R \geq 0$ such that $f(z) = \sum_{n\geq 0} z^n a_n$ $(z \in \mathbb{C})$ converges a.s. for $|z| < R$ and diverges a.s. for $|z| > R$. The function $z \to f(z)$ is a.s. analytic for $|z| < R$ and the convergence is uniform on every compact $\subset \mathcal{B}(0, R)$.
2. The radius of convergence $R$ is given by

$$\frac{1}{R} = \sup\left\{A > 0 \Big| \sum_{n\geq 0} \mathbb{P}[|a_n| \geq A^n] = \infty\right\}.$$

PROOF.

1. By Hadamard's criterion, the radius of convergence $R$ of a series $\sum_{n\geq 0} a_n z^n$ is given by

$$\frac{1}{R} = \limsup_n \sqrt[n]{|a_n|}.$$

Kolmogorov's (0–1) law states that $\mathbb{P}[\limsup \sqrt[n]{|a_n|} \geq A] = 0$ or 1 since this is a tail event. Hence, $1/R = \sup\{A > 0 | \mathbb{P}[\limsup \sqrt[n]{|a_n|} \geq A] = 1\}$ if this set is nonempty, or 0 else.
2. By the Borel–Cantelli lemma,

$$\mathbb{P}[\limsup \sqrt[n]{|a_n|} \geq A] = \mathbb{P}\left[\bigcap_{n\geq 0} \bigcup_{k\geq n} \{|a_k| \geq A^k\}\right] = 1$$

if and only if $\sum_{n\geq 0} \mathbb{P}[|a_n| \geq A^n] = \infty$. □

COROLLARY 1.4. *Suppose $f(z) = \sum_{n\geq 0} \lambda_n z^n \xi_n$, where the $(\xi_n)$ are independent standard Gaussian variables. Let $1/R := \limsup_{n\to\infty} \sqrt[n]{|\lambda_n|}$. Then the series converges uniformly on every compact $\subset \mathcal{B}(0, R)$ to a holomorphic function.*

---

[1]One should actually set $\tilde{\Gamma}'^+(\zeta) = \Gamma'^+(z)\frac{dz}{d\zeta} = \frac{2i}{(1-\zeta)^2}\Gamma'(\mathcal{C}^{-1}(\zeta))$ since $\Gamma'^+$ is to be interpreted as a derivative process. But this is not important for what follows.

PROOF. Supposing $A > 1/R$, then $\mathbb{P}[|\lambda_n \xi_n| \geq A^n] \leq C(\frac{e^{-A^{2n}/2\lambda_n^2}}{A^n/\lambda_n})$ goes rapidly to 0 and the series converges, while it is clear that $\sum_n \mathbb{P}[|\lambda_n \xi_n| \geq A^n] = \infty$ if $A < 1/R$. □

It follows immediately from the preceding corollary and (1.3) that $\Gamma'^+_z$ is well defined and analytic for $z \in \Pi^+$.

The next technical lemma is crucial to establish the connection with FBM and regularity properties.

LEMMA 1.5. *Let $\gamma, \gamma' : (0,1) \to \Pi^+$ be continuous paths with endpoints $\gamma(0) = a_1 + ib_1, \gamma(1) = a_2 + ib_2, \gamma'(0) = a'_1 + ib'_1, \gamma'(1) = a'_2 + ib'_2$ in the closed upper half-plane $\bar{\Pi}^+$, that is, $b_1, b_2, b'_1, b'_2 \geq 0$. Then, for every $\alpha \in (0,1)$, $\alpha \neq \frac{1}{2}$:*

1. *The double integral*

$$I = \int_\gamma du \int_{\bar{\gamma}'} d\bar{v} (-i(u-\bar{v}))^{2\alpha-2}$$

*is well defined. It depends only on the homology class of the paths $\gamma, \gamma'$, or, in other words, on the endpoints $\gamma(0), \gamma'(0), \gamma(1), \gamma'(1)$. Furthermore, it is invariant under real translations $\gamma \to \gamma + a, \gamma' \to \gamma' + a$ ($a \in \mathbb{R}$).*

2. *Suppose that $a_1 = a'_1 = 0$ and take $b_1, b_2, b'_1, b'_2 \to 0$. Then*

$$I_{\gamma, \gamma'} \to \frac{e^{i\pi\alpha \operatorname{sgn}(a'_2 - a_2)}|a'_2 - a_2|^{2\alpha} - e^{-i\pi\alpha \operatorname{sgn}(a_2)}|a_2|^{2\alpha} - e^{i\pi\alpha \operatorname{sgn}(a'_2)}|a'_2|^{2\alpha}}{2\alpha(2\alpha - 1)}.$$

*Hence,*

$$I_{\gamma, \gamma'} + \overline{I_{\gamma, \gamma'}} = 2 \operatorname{Re} I_{\gamma, \gamma'} \to \frac{2\cos\pi\alpha}{\alpha(1-2\alpha)} \cdot \frac{|a_2|^{2\alpha} + |a'_2|^{2\alpha} - |a'_2 - a_2|^{2\alpha}}{2}.$$

3. *Suppose $\gamma = \gamma'$. Then there exists a constant $C > 0$ such that*

$$\left| \int_\gamma du \int_{\bar{\gamma}} d\bar{v} (-i(u-\bar{v}))^{2\alpha-2} \right| \leq C |\gamma(1) - \gamma(0)|^{2\alpha}.$$

PROOF. Point 1 follows from Cauchy's theorem and the fact that

$$|(-i(u-\bar{v}))^{2\alpha-2}| \leq (\operatorname{Im} u + \operatorname{Im} v)^{2\alpha-2}$$

and

$$\int_0^1 dx \int_0^1 dy (x+y)^{2\alpha-2} < \infty.$$

As for point 2, an explicit computation yields (if $a_1 = a'_1 = 0$)

$$\int_\gamma \int_{\bar\gamma'} (-i(u-\bar v))^{2\alpha-2} \, du \, d\bar v$$

$$= ((b_2 + b'_2 - i(a_2 - a'_2))^{2\alpha} - (b_2 + b'_1 - ia_2)^{2\alpha}$$
$$- (b_1 + b'_2 + ia'_2)^{2\alpha} + (b_1 + b'_1)^{2\alpha})(2\alpha(2\alpha-1))^{-1}.$$

When $b_1, b_2, b'_1, b'_2 \to 0$, this expression goes to

$$\frac{e^{i\pi\alpha\,\mathrm{sgn}(a'_2 - a_2)}|a'_2 - a_2|^{2\alpha} - e^{-i\pi\alpha\,\mathrm{sgn}(a_2)}|a_2|^{2\alpha} - e^{i\pi\alpha\,\mathrm{sgn}(a'_2)}|a'_2|^{2\alpha}}{2\alpha(2\alpha-1)},$$

hence, also the result for $\mathrm{Re}\, I_{\gamma,\gamma'}$.

So let us now set about to show point 3. We split the proof into several parts.

1. Suppose $\gamma(0) = 0$ and $\gamma(1) = z$, $\mathrm{Im}\, z \geq 0$. Then one finds

(1.4)
$$\int_\gamma \int_{\bar\gamma} du \, d\bar v (-i(u-\bar v))^{2\alpha-2} = (2\,\mathrm{Im}\, z)^{2\alpha} - 2\,\mathrm{Re}(-iz)^{2\alpha}$$
$$= (2r\sin\theta)^{2\alpha} - 2\,\mathrm{Re}\, r^{2\alpha} e^{2i\alpha(\theta-\pi/2)},$$

if $z = re^{i\theta}, \theta \in [0, \pi]$, hence, there exists a constant $C > 0$ such that

$$\left| \int_\gamma \int_{\bar\gamma} du \, d\bar v (-i(u-\bar v))^{2\alpha-2} \right| \leq C|z|^{2\alpha}.$$

By invariance under real translations $z \to z + a$, $a \in \mathbb{R}$, a similar estimate holds for $\gamma(0)$ real.

2. Suppose $\frac{1}{2}b_2 \leq b_1 \leq b_2$. Set $u = x_1 + iy_1, v = x_2 + iy_2$. Let us compute the above double integral $I_{\gamma,\gamma}$ by choosing $\gamma = \gamma_1 \cup \gamma_2 := [a_1 + ib_1, a_2 + ib_1] \cup [a_2 + ib_1, a_2 + ib_2]$ as a horizontal line followed by a vertical line.

The double integral on the vertical line $\int\int_{\gamma_2 \times \bar\gamma_2} du \, d\bar v (-i(u-\bar v))^{2\alpha-2}$ is less than

(1.5)
$$\int_{b_1}^{b_2} \int_{b_1}^{b_2} dy_1 \, dy_2 (2b_1)^{2\alpha-2}$$
$$\leq C b_1^{2\alpha-2}(b_2 - b_1)^2 = \left( \frac{|b_2 - b_1|}{b_1} \right)^{2-2\alpha} |b_2 - b_1|^{2\alpha}$$
$$\leq |b_2 - b_1|^{2\alpha}.$$

The mixed integrals

$$I_{ij} := \int_{\gamma_i} \int_{\bar\gamma_j} du \, d\bar v (-i(u-\bar v))^{2\alpha-2}, \qquad i \neq j,$$

satisfy

$$|I_{ij}| \leq C \int_0^{|a_2-a_1|} dx \int_0^{b_2-b_1} dy (2b_1 + x + y)^{2\alpha-2}$$
$$\leq C'|a_1 - a_2|(b_2 - b_1)b_1^{2\alpha-2}, \tag{1.6}$$

hence, also

$$|I_{ij}| \leq C'|a_1 - a_2|^{2\alpha} \left(\frac{b_1}{a_1 - a_2}\right)^{2\alpha-1}. \tag{1.7}$$

On the other hand,

$$\int_0^{|a_2-a_1|} dx \int_0^{b_2-b_1} dy (2b_1 + x + y)^{2\alpha-2}$$
$$= \frac{1}{2\alpha(2\alpha-1)}[(b_1 + b_2 + |a_2 - a_1|)^{2\alpha} - (b_1 + b_2)^{2\alpha} \tag{1.8}$$
$$- (2b_1 + |a_2 - a_1|)^{2\alpha} + (2b_1)^{2\alpha}].$$

We must now distinguish 4 cases according to the sign of $\alpha - \frac{1}{2}$ and to the relative order of magnitude of $b_1$ and $|a_1 - a_2|$.

Suppose first $\alpha > \frac{1}{2}$: then $|I_{ij}| \leq C'|a_1 - a_2|^{2\alpha}$ by (1.7) if $b_1 \leq |a_1 - a_2|$, and $|I_{ij}| \leq C'|a_1 - a_2|^{2\alpha-1}(b_2 - b_1)$ if $b_1 \geq |a_1 - a_2|$ by (1.6). Hölder's inequality

$$|xy| \leq \frac{|x|^p}{p} + \frac{|y|^q}{q}, \qquad \frac{1}{p} + \frac{1}{q} = 1$$

applied to $x = |a_1 - a_2|^{(2\alpha-1)/\alpha}$, $y = (b_2 - b_1)^{1/\alpha}$, $p = \frac{2\alpha}{2\alpha-1}$ yields in the latter case $|I_{ij}| \leq C''(|a_1 - a_2|^2 + |b_1 - b_2|^2)^\alpha$.

If now $\alpha < \frac{1}{2}$, then $|I_{ij}| \leq C'|a_1 - a_2|^{2\alpha}$ by (1.7) if $\frac{b_1}{|a_1-a_2|}$ is bounded below by a positive constant, while (1.8) yields for $|b_1| \ll |a_1 - a_2|$

$$|I_{ij}| \leq C \sup(b_1^{2\alpha}, |a_1 - a_2|^{2\alpha-1}(b_2 - b_1)) \ll |a_1 - a_2|^{2\alpha}.$$

Finally, the integral on the horizontal line

$$I_{11} = \int_{\gamma_1} \int_{\bar{\gamma}_1} du\, d\bar{v}(-i(u-\bar{v}))^{2\alpha-2} = \int_{a_1}^{a_2} \int_{a_1}^{a_2} dx_1\, dx_2 (2b_1 - i(x_1 - x_2))^{2\alpha-2}$$

can be evaluated explicitly, namely,

$$I_{11} = \frac{2(2b_1)^{2\alpha} - (2b_1 - i(a_1 - a_2))^{2\alpha} - (2b_1 + i(a_1 - a_2))^{2\alpha}}{2\alpha(2\alpha - 1)}.$$

An expansion at order two by using Taylor's formula gives $I_{11} \leq Cb_1^{2\alpha-2}|a_1 - a_2|^2 \leq |a_1 - a_2|^{2\alpha}$ if $b_1 \gg |a_1 - a_2|$. If $b_1$ is of the same order as $a_1 - a_2$ or smaller, the same estimate comes out trivially.

3. Suppose now that $b_2 > 2b_1$. Choose the following contour of integration: set $\gamma = \gamma_1 \amalg \gamma_2$, with $\gamma_1 = \{a_1 + ib_1(1-t) | t \in [0,1]\}$ and $\gamma_2 = \{(1-t)a_1 + t(a_2 + ib_2) | t \in [0,1]\}$. Then

$$\int_\gamma \int_{\bar\gamma} du\, d\bar v (-i(u-\bar v))^{2\alpha-2} = I_{11} + I_{12} + I_{21} + I_{22},$$

where $I_{jk} = \int_{\gamma_j} \int_{\gamma_k} du\, d\bar v (-i(u-\bar v))^{2\alpha-2}$. The integral $I_{11}$ (resp. $I_{22}$) is less that a constant times $b_1^{2\alpha}$ [resp. $((a_1-a_2)^2 + b_2^2)^\alpha$] by point 1. As for the mixed terms $I_{ij}, i \neq j$, an explicit computation yields, for instance,

$$I_{12} = \frac{1}{2\alpha(2\alpha-1)}[(b_1+b_2-i(a_1-a_2))^{2\alpha} - b_1^{2\alpha} - (b_2 - i(a_1-a_2))^{2\alpha}],$$

which is $\leq C(|a_1-a_2|^2 + |b_1-b_2|^2)^\alpha$. By putting everything together, one gets

$$\left| \int_\gamma \int_{\bar\gamma} du\, d\bar v (-i(u-\bar v))^{2\alpha-2} \right| \leq C'((a_1-a_2)^2 + (b_1-b_2)^2)^\alpha. \qquad \square$$

DEFINITION 1.6. Define for any $z \in \bar\Pi^+$

$$\Gamma_z^+ = \int_\gamma du\, \Gamma_u'^+,$$

where $\gamma : (0,1) \to \Pi^+$ is any continuous path with endpoints $\gamma(0) = 0, \gamma(1) = z$ and, similarly, $\Gamma_{\bar z}^- = \overline{(\Gamma_z^+)}$.

Lemma 1.5 shows that $\Gamma_z^\pm$ has a well defined covariance on $\bar\Pi^\pm$. One has thus obtained the $\Gamma^\pm$-*processes*. These are hermitian Gaussian processes defined on $\bar\Pi^\pm$. Restricted to $\Pi^\pm$, they have analytic paths.

COROLLARY 1.7.

1. *The real-time real-valued process* $(\Gamma_t, t \geq 0)$ *defined as*

$$\Gamma_t = 2\operatorname{Re}\Gamma_t^+ = \Gamma_t^+ + \Gamma_t^-$$

*has the same law as* $(B_t, t \in \mathbb{R})$.
2. *Let* $t \in \mathbb{R}$. *Then, with probability* 1,

$$2\operatorname{Re}\Gamma_z^+ = \Gamma_z^+ + \Gamma_{\bar z}^- \to \Gamma_t$$

*as* $z \in \Pi^+$ *goes to* $t$ *along any path in the closed upper half-plane. More precisely,*

$$\mathbb{E}[|\operatorname{Re}\Gamma_z^+ - \Gamma_t|^2] \leq C|z-t|^{2\alpha}.$$

*Hence, the real-valued process* $z \to \Gamma_z^+$, $z \in \Pi^+$ *has an* $(\alpha - \varepsilon)$-*Hölder continuous version for every* $\varepsilon > 0$.

REMARK. In other words, FBM appears as the boundary value of the analytic process defined as $\Gamma_z^+$ on $\Pi^+$ and $\Gamma_{\bar{z}}^-$ on $\Pi^-$. We shall come back to this in Section 3.

PROOF. It is straightforward from Lemma 1.5 by using Kolmogorov's lemma for Gaussian processes. □

We shall see in Section 3 that the convergence of $\operatorname{Re}\Gamma_z^+$ to FBM as $\operatorname{Im} z \to 0$ is a.s. uniform on every compact.

But let us first come back to the series $\Gamma_z'^+ = \sum_{k\geq 0} f_k(z)\xi_k^+$ ($z \in \Pi^+$) and show how FBM appears as the limit of a series defined directly on $\mathbb{R}$.

**2. On a series decomposition of FBM.** Let us recall the definition of $f_k(z), z \in \Pi^+$ ($k \geq 0$) from Lemma 1.1 in Section 1:

$$f_k(z) = 2^{\alpha-1}\sqrt{\frac{\alpha(1-2\alpha)}{2\cos\pi\alpha}}\sqrt{\frac{(2-2\alpha)_k}{k!}}\left(\frac{z+i}{2i}\right)^{2\alpha-2}\left(\frac{z-i}{z+i}\right)^k.$$

This function extends analytically to the half-plane $\operatorname{Im} z > -1$; in any case, it is well defined on the real axis. By Stirling's formula,

$$\sqrt{\frac{(2-2\alpha)_k}{k!}} = \sqrt{\frac{1}{\Gamma(2-2\alpha)}}\sqrt{\frac{\Gamma(2-2\alpha+k)}{\Gamma(1+k)}} = O_{k\to\infty}(k^{1/2-\alpha}).$$

Recall also that $|\frac{z-i}{z+i}| = 1$ for $z \in \mathbb{R}$. Hence, the series $\sum_{k\geq 0} |f_k(t)|^2$ diverges on the real axis as $\sum_{k\geq 1} k^{1-2\alpha}$, which is of course not at all surprising [otherwise the random series $\sum_{k\geq 0} f_k(t)\xi_k^+ + \sum_{k\geq 0} \overline{f_k(t)}\xi_k^-$ would have been a good candidate for a would-be derivative of FBM!].

Set

$$F_k(z) := \int_0^z f_k(s)\,ds, \qquad z \in \bar{\Pi}^+.$$

LEMMA 2.1.

1. Let $C_\alpha(s,t) = \sum_{k\geq 0} F_k(s)\overline{F_k(t)}$, $s,t \in \mathbb{R}$. Then this series converges in absolute value.
2. One has $C_\alpha(s,t) = \frac{e^{-i\pi\alpha\operatorname{sgn}(s)}|s|^{2\alpha}+e^{i\pi\alpha\operatorname{sgn}(t)}|t|^{2\alpha}-e^{i\pi\alpha\operatorname{sgn}(t-s)}|t-s|^{2\alpha}}{4\cos\pi\alpha}$.
3. Let $B_t := \sum_{k\geq 0} F_k(t)\xi_k^+ + \sum_{k\geq 0} \overline{F_k(t)}\xi_k^-$. Then $(B_t, t \in \mathbb{R})$ is a two-sided FBM with Hurst index $\alpha$.

PROOF. Let $z$ stay in a fixed compact set $K \subset \bar{\Pi}^+$. Then

$$\int_0^z du\left(\frac{u-i}{u+i}\right)^k\left(\frac{u+i}{2i}\right)^{2\alpha-2}$$

$$= \frac{1}{k+1} \int_0^z du \frac{d}{du}\left[\left(\frac{u-i}{u+i}\right)^{k+1}\right] \frac{(u+i)^{2\alpha}}{(2i)^{2\alpha-1}}$$

$$= \frac{1}{k+1}(2i)^{1-2\alpha}\left\{\left(\frac{z-i}{z+i}\right)^{k+1}(z+i)^{2\alpha} - (-1)^{k+1}e^{i\pi\alpha}\right.$$

$$\left. - 2\alpha \int_0^z du \left(\frac{u-i}{u+i}\right)^{k+1}(u+i)^{2\alpha-1}\right\}.$$

Hence,

$$|F_k(z)| \leq C\sqrt{\frac{(2-2\alpha)_k}{k!}} \frac{1}{k+1} \leq (1+k)^{-1/2-\alpha},$$

which guarantees that the series $C_\alpha(z,\bar{w}) := \sum_{k\geq 0} F_k(z)\overline{F_k(w)}$ is normally convergent in $K$. Consequently,

$$C_\alpha(s,t) = \lim_{z\to s, \bar{w}\to t} C_\alpha(z,\bar{w}),$$

if $z \to s, w \to t$ while staying inside $\Pi^+$. But (see Lemma 1.2) the series $\sum f_k(z)\overline{f_k(w)}$ converges in absolute value if $z, w \in \Pi^+$, hence (by commuting the sum with the integral),

$$C_\alpha(z,\bar{w}) = \frac{\alpha(1-2\alpha)}{2\cos\pi\alpha} \int_0^z du \int_0^{\bar{w}} d\bar{v}(-i(u-\bar{v}))^{2\alpha-2}.$$

Now one may conclude by using Lemma 1.5(2). □

By standard arguments using martingale inequlities and Fernique's lemma (see, e.g., Theorem 3.3.2 in [2]), the Karhunen–Loève type series

(2.1) $$t \to \sum_{k\geq 0} F_k(t)\xi_k^+ + \sum_{k\geq 0} \overline{F_k(t)}\xi_k^-$$

converges uniformly with probability one to $t \to B_t$ on every compact set $K \subset \mathbb{R}$.

Let us now study the rate of convergence of the above series. It may be computed by using standard entropy arguments that we reproduce from [2]. Set $R_N(t) := \sum_{k\geq N} F_k(t)\xi_k^+ + \sum_{k\geq N} \overline{F_k(t)}\xi_k^-$. Then the standard deviation semi-metric

$$d_N(s,t) := \sqrt{\mathbb{E}[(R_N(t)-R_N(s))^2]}$$

defines a natural metric on $\mathbb{R}$ from the point of view of the Gaussian process $t \to R_N(t)$. We shall be using the following

PROPOSITION 2.2 (See [20], page 101). *Let $X_t$ be a centered Gaussian process on $K \subset \mathbb{R}$ compact, and $d(s,t) := \sqrt{\mathbb{E}[|X_t - X_s|^2]}$ its associated (semi-)metric. Then*

$$\mathbb{E}\left[\sup_{t \in K} |X_t|\right] \leq C \int_0^{+\infty} \sqrt{\ln N(d, K, \varepsilon)}\, d\varepsilon,$$

*where $C$ is a universal constant, and $N(d, K, \varepsilon)$ is the smallest number of balls of $d$-radius $\varepsilon$ covering $K$.*

The logarithm of $N(d, K, \varepsilon)$ is usually called the entropy of the process $X$ on $K$.

Let $s, t$ remain in a fixed compact $K \subset \mathbb{R}$. Then $|R_N(t) - R_N(s)|$ may be estimated either by remarking that

$$\mathbb{E}[|R_N(t) - R_N(s)|^2] = 2 \sum_{k \geq N} |F_k(t) - F_k(s)|^2 \leq 2 \sum_{k \geq 0} |F_k(t) - F_k(s)|^2$$
(2.2)
$$= \mathbb{E}[|B_t - B_s|^2] = |t - s|^{2\alpha}$$

or by using an upper bound appearing in the proof of Lemma 2.1, namely,

$$\mathbb{E}[|R_N(t) - R_N(s)|^2] \leq 8 \sum_{k \geq N} \sup(|F_k(t)|^2, |F_k(s)|^2)$$
(2.3)
$$\leq C \sum_{k \geq N} k^{-1-2\alpha} = C' N^{-2\alpha}.$$

Hence (the following trick is borrowed from a similar proof in [7]), for every couple $(p, q)$ of positive numbers,

(2.4) $$\mathbb{E}[|R_N(t) - R_N(s)|^2]^{p+q} \leq C'' N^{-2p\alpha} |t - s|^{2q\alpha}.$$

Consequently, $\ln N(d, K, \varepsilon) \sim \ln(\varepsilon^{-(p+q)/q\alpha} N^{-p/q})$ for $\varepsilon \lesssim N^{-p\alpha/(p+q)}$, and $\ln N(d_B, K, \varepsilon) = 0$ for $\varepsilon \gtrsim N^{-p\alpha/(p+q)}$, leading to

$$\mathbb{E}\left[\sup_{t \in K} R_N(t)\right] \leq C \sqrt{\frac{p+q}{q\alpha}} N^{-p\alpha/p+q} \int_0^1 \sqrt{\ln(1/x)}\, dx \leq C' \sqrt{\frac{p+q}{q\alpha}} N^{-p\alpha/p+q}.$$

The best estimate for $\sup_{t \in K} |R_N(t)|$ is obtained in the limit $p \to \infty$ (with $q$ fixed), namely

$$\mathbb{E}\left[\sup_{t \in K} |R_N(t)|\right] \leq C_\varepsilon N^{-(\alpha - \varepsilon)}$$

for every $\varepsilon > 0$, with a constant $C_\varepsilon$ that diverges when $\varepsilon \to 0$.

Summarizing, we have proved the following:

THEOREM 2.3. *Let $K \subset \mathbb{R}$ be a compact. Then, for every $\varepsilon > 0$,*

$$N^{-\varepsilon+\alpha}\mathbb{E}\left[\sup_K \left| B_t - \sum_{k \leq N} F_k(t)\xi_k^+ - \sum_{k \leq N} \overline{F_k(t)}\xi_k^- \right|\right] \stackrel{N \to \infty}{\to} 0.$$

Another unrelated series expansion of the fractional Brownian motion (using zeros of Bessel functions) has been obtained by Dzhaparidze and Van Zanten in 2004 (see [7]). The convergence appeared to be *optimal* in the sense of a paper by Kühn and Linde (see [11]): they proved that no series converging uniformly on every compact to $\text{FBM}_\alpha$ could have a remainder of order less than $N^{-\alpha}\sqrt{\log N}$. Hence, the convergence is considered to be optimal if the remainder is $\leq C_\varepsilon N^{(-\alpha-\varepsilon)}$ for every $\varepsilon > 0$. We refer to these two papers for details. Here also, as we have just shown, the convergence is optimal.

**3. First approximation of FBM.** From the results of Section 1, it is natural to approximate FBM by the following process:

DEFINITION 3.1. For any $\varepsilon > 0$, let $(\Gamma(\varepsilon)_t, t \in \mathbb{R})$ be the Gaussian process defined by

$$\Gamma(\varepsilon)_t = \operatorname{Re} \Gamma^+(t + i\varepsilon) = \Gamma^+(t + i\varepsilon) + \Gamma^-(t - i\varepsilon).$$

By Lemma 1.5, $\Gamma(\varepsilon)_t$ has a (Gaussian) limit $\Gamma_t$ in $L^2$ as $\varepsilon \to 0$, and the limit process $t \to \Gamma_t$ has the same law as FBM. In other words, viewing $\Gamma^+$ and $\Gamma^-$ as the expression on $\Pi^+$ and on $\Pi^-$ of the same analytic process $\Gamma(z)$ defined on $\mathbb{C} \setminus \mathbb{R}$, FBM may be seen as the *boundary value* of $\Gamma$, namely,

(3.1) $$\Gamma_t := \lim_{\varepsilon \to 0}(\Gamma(t + i\varepsilon) + \Gamma(t - i\varepsilon)).$$

More precisely, still by Lemma 1.5, there exists a constant $C > 0$ such that, for every $t \in \mathbb{R}$,

$$\mathbb{E}[|\Gamma(\varepsilon)_t - \Gamma_t|^2] \leq C\varepsilon^{2\alpha}.$$

3.1. *Uniform approximation of FBM by $\Gamma(\varepsilon)$.* We shall give in this paragraph, by a method similar to that of the preceding section, an estimate of the sup-norm of the difference $\Gamma - \Gamma(\varepsilon)$.

THEOREM 3.2. *Let $K \subset \mathbb{R}$ be a compact. Then, for every $\eta > 0$,*

$$\varepsilon^{\eta-\alpha}\mathbb{E}\left[\sup_K |\Gamma_t - \Gamma(\varepsilon)_t|\right] \stackrel{\varepsilon \to 0}{\to} 0.$$

PROOF. The proof goes as in Theorem 2.3. Namely, the estimate (2.2) may be replaced by

$$\mathbb{E}[|(\Gamma_t - \Gamma(\varepsilon)_t) - (\Gamma_s - \Gamma(\varepsilon)_s)|^2]$$
$$\leq 4\sup(\mathbb{E}[|\Gamma_t - \Gamma_s|^2], \mathbb{E}[|\Gamma(\varepsilon)_t - \Gamma(\varepsilon)_s|^2])$$
$$\leq C|t-s|^{2\alpha}$$

(see Lemma 1.5) and the estimate (2.3) by

$$\mathbb{E}[|(\Gamma_t - \Gamma(\varepsilon)_t) - (\Gamma_s - \Gamma(\varepsilon)_s)|^2]$$
$$\leq 4\sup(\mathbb{E}[|\Gamma_t - \Gamma(\varepsilon)_t|^2], \mathbb{E}[|\Gamma_s - \Gamma(\varepsilon)_s|^2])$$
$$\leq C'\varepsilon^{2\alpha}$$

(also by Lemma 1.5). $\square$

3.2. *Convergence of analytic iterated integrals.* Obtaining good estimates of iterated integrals is the key point in order to apply rough path theory. We claim that the *analytic* iterated integrals, as defined below, converge in $L^2(\mathbb{P})$ for *every* Hurst index $\alpha > 0$.

DEFINITION 3.3. Let $s,t > 0$ and $f_1, \ldots, f_n$ ($n \geq 1$) be $n$ analytic functions defined on a neighborhood of the closed strip $\Omega = \{z = x+iy | 0 \leq x \leq t, 0 \leq y \leq s\}$. Then an *analytic iterated integral* is an integral of the form

$$\int_0^t dx_1 f(x_1 + i\varepsilon_1) \Gamma'(\varepsilon_1)_{x_1}^{+,(i_1)} \int_0^{x_1} dx_2$$
(3.2)
$$\cdots \int_0^{x_{n-1}} dx_n f(x_n + i\varepsilon_n) \Gamma'(\varepsilon_n)_{x_n}^{+,(i_n)},$$

where $\varepsilon_1, \ldots, \varepsilon_n > 0$ and $1 \leq i_1, \ldots, i_n \leq d$ (recall $d$ is the number of components of FBM).

THEOREM 3.4. *Let $s,t > 0$ and $f_1, \ldots, f_n$ (resp. $g_1, \ldots, g_n$) be $n$ analytic (resp. anti-analytic) functions defined on a neighborhood of the closed strip $\Omega = \{z = x+iy | 0 \leq x \leq t, 0 \leq y \leq s\}$ (resp. $\bar{\Omega}$). Then*

$$\mathcal{V}(\boldsymbol{\varepsilon}, \boldsymbol{\eta}) := \mathbb{E}\bigg[\bigg(\int_0^t dx_1 f_1(x_1 + i\varepsilon_1) \Gamma'(\varepsilon_1)_{x_1}^{+,(i_1)} \int_0^{x_1} dx_2$$
$$\cdots \int_0^{x_{n-1}} dx_n f_n(x_n + i\varepsilon_n) \Gamma'(\varepsilon_n)_{x_n}^{+,(i_n)}\bigg)$$
(3.3)
$$\times \bigg(\int_0^t dy_1 g_1(y_1 - i\varepsilon_1) \Gamma'(\eta_1)_{y_1}^{-,(j_1)} \int_0^{y_1} dy_2$$
$$\cdots \int_0^{y_{n-1}} dx_n g_n(y_n - i\eta_n) \Gamma'(\eta_n)_{y_n}^{-,(j_n)}\bigg)\bigg]$$

admits a limit $\mathcal{V}(0,0)$ when $\boldsymbol{\varepsilon} = (\varepsilon_1, \ldots, \varepsilon_n), \boldsymbol{\eta} = (\eta_1, \ldots, \eta_n)$ go to zero, which is bounded by

$$|\mathcal{V}(0,0)| \leq C \prod_{i=1}^{n} \sup_{z \in \Omega} |f_i(z)| \prod_{i=1}^{n} \sup_{\bar{w} \in \bar{\Omega}} |g_i(\bar{w})|$$

$$\times \max(s^{2\alpha}, ts^{2\alpha-1}, t^{2\alpha-2}s^2, t^2 s^{2\alpha-2})^n \qquad (3.4)$$

for a universal constant $C$.

PROOF. Supposing $i_1 \neq i_2 \neq \cdots \neq i_n$, one may set $j_1 = \sigma(i_1), \ldots, j_n = \sigma(i_n)$ for a certain permutation $\sigma$ [otherwise $\mathcal{V}(\varepsilon, \eta) = 0$]. Let $s' = s - \sup_k \varepsilon_k - \sup_k \eta_k$. Then

$$\mathcal{V}(\varepsilon, \eta) = \left(\frac{\alpha(1-2\alpha)}{2\cos \pi\alpha}\right)^n \left(\int_0^t dx_1 \int_0^{x_1} dx_2 \cdots \int_0^{x_{n-1}} dx_n\right)$$

$$\times \left(\int_0^t dy_1 \int_0^{y_1} dy_2 \cdots \int_0^{y_{n-1}} dy_n\right) \qquad (3.5)$$

$$\times \left(\prod_k f(x_k + i\varepsilon_k)\right)\left(\prod_k g(y_k - i\eta_k)\right)$$

$$\times \prod_k (-i(x_k - y_{\sigma(k)}) + \varepsilon_k + \eta_{\sigma(k)})^{2\alpha-2}$$

may be extended as a holomorphic function of $x_1, \ldots, x_n$ on the strip $\{z = x + iy | 0 \leq x \leq t, 0 \leq y \leq s - \sup_k \varepsilon_k\}$, respectively, as an anti-holomorphic function of $y_1, \ldots, y_n$ on the strip $\{\bar{z} = x - iy | 0 \leq x \leq t, 0 \leq y \leq s - \sup_k \eta_k\}$. Set

$$(3.6) \qquad \gamma_0 = [0, is'] \cup [is', t + is'] \cup [t + is', t], \qquad \gamma_\varepsilon = \gamma_0 + i\varepsilon$$

and, if $z \in \gamma_\varepsilon$, let $\gamma_\varepsilon(z)$ be the section of the path $\gamma_\varepsilon$ comprised between 0 and $z$. Then, by a deformation of contour, one gets

$$\mathcal{V}(\varepsilon, \eta) = \left(\frac{\alpha(1-2\alpha)}{2\cos \pi\alpha}\right)^n$$

$$\times \left(\int_{\gamma_{\varepsilon_1}} dz_1 \int_{\gamma_{\varepsilon_1}(z_1) + i(\varepsilon_2 - \varepsilon_1)} dz_2 \cdots \int_{\gamma_{\varepsilon_{n-1}}(z_{n-1}) + i(\varepsilon_n - \varepsilon_{n-1})} dz_n\right)$$

$$(3.7) \quad \times \left(\int_{\overline{\gamma_{\eta_1}}} d\bar{w}_1 \int_{\overline{\gamma_{\eta_1}(w_1) + i(\eta_2 - \eta_1)}} d\bar{w}_2 \cdots \int_{\overline{\gamma_{\eta_{n-1}}(w_{n-1}) + i(\eta_n - \eta_{n-1})}} d\bar{w}_n\right)$$

$$\times \left(\prod_k f(z_k)\right)\left(\prod_k g(\bar{w}_k)\right) \prod_k (-i(z_k - \bar{w}_{\sigma(k)}))^{2\alpha-2}.$$

Let $M := \prod_k \sup_{z \in \Omega} |f_k(z)| \prod_k \sup_{\bar{w} \in \bar{\Omega}} |g_k(\bar{w})|$. Then

$$
\begin{aligned}
|\mathcal{V}(\varepsilon, \eta)| &\leq C \cdot M \left( \int_{\gamma_0} |dz_1| \cdots \int_{\gamma_0} |dz_n| \right) \left( \int_{\bar{\gamma}_0} |d\bar{w}_1| \cdots \int_{\bar{\gamma}_0} |d\bar{w}_n| \right) \\
&\quad \times \prod_k |(-i(z_k - \bar{w}_{\sigma(k)}))^{2\alpha-2}|.
\end{aligned}
\tag{3.8}
$$

The different components decouple completely and one is left with the problem of estimating

$$
I := \int_{\gamma_0} \int_{\bar{\gamma}_0} |dz||d\bar{w}||(-i(z-\bar{w}))^{2\alpha-2}|.
\tag{3.9}
$$

Splitting $\gamma_0$ and $\bar{\gamma}_0$ into their three linear pieces, one sees that

$$
\begin{aligned}
\int_{[0,is]} |dz| \int_{[0,-is]} |d\bar{w}||(-i(z-\bar{w}))^{2\alpha-2}| &= \int_0^s dx \int_0^s dy (x+y)^{2\alpha-2} \\
&= \frac{2^{2\alpha} - 2}{2\alpha(2\alpha - 1)} s^{2\alpha};
\end{aligned}
\tag{3.10}
$$

$$
\int_{[is,is+t]} |dz| \int_{[0,-is]} |d\bar{w}||(-i(z-\bar{w}))^{2\alpha-2}| \leq t s^{2\alpha-1};
\tag{3.11}
$$

$$
\int_{[0,is]} |dz| \int_{[-is+t,t]} |d\bar{w}||(-i(z-\bar{w}))^{2\alpha-2}| \leq t^{2\alpha-2} s^2;
\tag{3.12}
$$

$$
\begin{aligned}
&\int_{[is,is+t]} |dz| \\
&\quad \times \int_{[-is,-is+t]} |d\bar{w}||(-i(z-\bar{w}))^{2\alpha-2}| \leq 2^{2\alpha-2} t^2 s^{2\alpha-2},
\end{aligned}
\tag{3.13}
$$

hence the result.

If now some of the indices $i_1, \ldots, i_n$ coincide, then $\mathcal{V}(\varepsilon, \eta)$ may be obtained as a finite linear combination of terms of the same type as (3.5), hence, the same estimates hold. (Easy) details are left to the reader. □

In particular, the second moment of any analytic iterated integral converges when $\varepsilon \to 0$. Hence the divergence of Lévy's area for FBM with Hurst index $\alpha \geq \frac{1}{4}$ stems from the mixed integral $\int_0^t dx_1 \Gamma'(\varepsilon_1)^{+,(i_1)} \int_0^{x_1} dx_2 \Gamma'(\varepsilon_2)^{-,(i_2)}$; see computation of the integral $\mathcal{V}_2$ in the proof of Theorem 4.4 below.

**4. Existence of the rough path limit for $\alpha > \frac{1}{4}$.** The purpose of this section is to prove that $\Gamma(\varepsilon)$ ($d$-dimensional $\varepsilon$-regularized FBM) has a rough path limits when $\varepsilon \to 0$ for every $d \geq 1$ and $\alpha > \frac{1}{4}$. Let us be a little more specific about what has to be done:

– for every $\alpha > 0$, $d_1(\Gamma(\varepsilon), \Gamma)_{s,t} \to 0$ with probability one (which is a simple consequence of Lemma 1.5);

– if $\alpha > \frac{1}{2}$, then no regularization is needed. Hence, we assume (for the proof of the rough path convergence) that $\alpha \in (\frac{1}{4}, \frac{1}{2})$. Yet the explicit computation (see Section 4.1) of the $\varepsilon$-regularized Levy area for FBM is valid for every Hurst index $\alpha$;

– if $\frac{1}{3} < \alpha < \frac{1}{2}$, one needs to prove that $\mathbf{\Gamma}(\varepsilon)$ converges in the $d_2$-distance for $q > 1/\alpha$. Since a simple integration by parts yields $\Gamma(\varepsilon)_{i,i}^2 = \frac{1}{2}(\mathbf{\Gamma}(\varepsilon)^{(i)})^2$, the only real problem is to prove the convergence of $\mathbf{\Gamma}(\varepsilon)_{i,j}^2$ for $i \neq j$. This we do in all details in Section 4.2; we shall actually prove that convergence in the $d_2$-distance holds for every $\alpha > \frac{1}{4}$;

– if $\frac{1}{4} < \alpha < \frac{1}{3}$, there remains then to prove the convergence of the third iterated integrals $\mathbf{\Gamma}(\varepsilon)_{i,j,k}^3$ in the $d_3$-distance. Since $\mathbf{\Gamma}(\varepsilon)_{i,i,i}^3 = \frac{1}{6}(\mathbf{\Gamma}(\varepsilon)^{(i)})^3$ and $\mathbf{\Gamma}(\varepsilon)_{i,i,j}^3$, $\mathbf{\Gamma}(\varepsilon)_{i,j,i}^3$, $\mathbf{\Gamma}(\varepsilon)_{j,i,i}^3$ ($i \neq j$) reduce easily (by integration by parts) to second-order iterated integrals of one of the following types,

$$\int_0^t (\Gamma(\varepsilon)_x^{(i)})^2 \Gamma'(\varepsilon)_x^{(j)}\, dx,$$

$$\Gamma(\varepsilon)_t^{(i)} \int_0^t \Gamma(\varepsilon)_x^{(i)} \Gamma'(\varepsilon)_x^{(j)}\, dx,$$

$$\Gamma(\varepsilon)_t^{(i)} \int_0^t \Gamma(\varepsilon)_x^{(j)} \Gamma'(\varepsilon)_x^{(i)}\, dx,$$

we shall only need to consider the case $i \neq j \neq k$.

4.1. *Computation of two-dimensional Lévy's area for $\alpha \in (0,1)$.* Generally speaking, given a two-component stochastic process $\mathbf{X}_t = (X_t^{(1)}, X_t^{(2)})$ with regular enough trajectories (say, for smooth trajectories), we call *Lévy's area*[2] of $\mathbf{X}$ between time $s$ and time $t$ the random quantity

$$\text{Area}_{s,t}(\mathbf{X}) = \int_s^t (X_u^{(2)} - X_s^{(2)})\, dX_u^{(1)}$$

(4.1)
$$= \int_s^t (X_u^{(2)} - X_s^{(2)}) \frac{d}{du} X_u^{(1)}\, du$$

$$= \int_s^t X_x'^{(1)}\, dx \int_s^x X_y'^{(2)}\, dy.$$

---

[2] The usual definition of Lévy's area is antisymmetric in $X^{(1)}$ and $X^{(2)}$, given by $\int_s^t X_u^{(2)}\, dX_u^{(1)} - X_u^{(1)}\, dX_u^{(2)}$. It differs from our definition in an inessential way for our purposes.

It is given geometrically as the area comprised between the curve

$$u \to \begin{pmatrix} x(u) \\ y(u) \end{pmatrix} := \begin{pmatrix} X_u^{(1)} \\ X_u^{(2)} \end{pmatrix}$$

and the union of two broken lines, respectively, horizontal and vertical, namely $[\binom{x(s)}{y(s)}, \binom{x(t)}{y(s)}] \cup [\binom{x(t)}{y(s)}, \binom{x(t)}{y(t)}]$.

The purpose of this paragraph is to give an explicit expression of this quantity for $\mathbf{X} = \mathbf{\Gamma}(\varepsilon)$, where $\mathbf{\Gamma}$ is a two-component FBM.

We shall be using a number of times the following integral representation of Gauss' hypergeometric function $_2F_1$ (see [1], (15.3.1)):

$$(4.2) \quad {}_2F_1(a,b,c;z) = \frac{\Gamma(c)}{\Gamma(b)\Gamma(c-b)} \int_0^1 t^{b-1}(1-t)^{c-b-1}(1-tz)^{-a}\,dt,$$

valid if $\operatorname{Re} c > \operatorname{Re} b > 0$ and $z \in \mathbb{C} \setminus [1; +\infty)$. Recall that $_2F_1(a,b,c;z)$ ($c \neq 0, -1, \ldots$) is defined around $z = 0$ by an infinite series with radius of convergence 1 and has an analytic extension to the cut plane $\mathbb{C} \setminus [1; +\infty[$. The *connection formulas* give $_2F_1(a,b,c;z)$ in terms of a linear combination of hypergeometric functions in the transformed argument $\phi(z)$, where $\phi$ is any projective transformation of the Riemann sphere preserving the set of singularities of the hypergeometric differential equation, namely, $\{0, 1, \infty\}$. They relate the behavior of the hypergeometric functions around 0 with their behavior around 1 and $\infty$. We reproduce here for the convenience of the reader two connection formulas, relating the behavior around 0 with the behavior around $\infty$, respectively, 1 (see [1], Section 15.3):

$$(4.3) \quad \begin{aligned} & {}_2F_1(a,b,c;z) \\ &= \frac{\Gamma(c)\Gamma(b-a)}{\Gamma(b)\Gamma(c-a)}(-z)^{-a}{}_2F_1\left(a, 1-c+a; 1-b+a; \frac{1}{z}\right) \\ &\quad + \frac{\Gamma(c)\Gamma(a-b)}{\Gamma(a)\Gamma(c-b)}(-z)^{-b}{}_2F_1\left(b, 1-c+b; 1-a+b; \frac{1}{z}\right), \end{aligned}$$
$$z \notin \mathbb{R}_+;$$

$$(4.4) \quad \begin{aligned} & {}_2F_1(a,b,c;z) \\ &= \frac{\Gamma(c)\Gamma(c-a-b)}{\Gamma(c-a)\Gamma(c-b)}{}_2F_1(a,b;a+b-c+1;1-z) \\ &\quad + (1-z)^{c-a-b}\frac{\Gamma(c)\Gamma(a+b-c)}{\Gamma(a)\Gamma(b)} \\ &\qquad \times {}_2F_1(c-a,c-b;c-a-b+1;1-z), \qquad z \notin [1; +\infty). \end{aligned}$$

Let us also recall that $_2F_1(a,b,c;z)$ is symmetric in the arguments $a,b$, constant (equal to 1) if $a=0$ or $b=0$, and that $\lim_{s\to 1, s\notin [1;+\infty[}{}_2F_1(a,b,c;s)$ exists and is finite if and only if $\text{Re}(c-a-b)>0$, in which case

$$(4.5) \quad {}_2F_1(a,b,c;1) = \frac{\Gamma(c)\Gamma(c-a-b)}{\Gamma(c-a)\Gamma(c-b)}, \qquad c\notin \mathbb{Z}_-, \text{Re}(c-a-b)>0.$$

LEMMA 4.1. *Let $\beta_1 \in \mathbb{C}, \beta_2 \in \mathbb{C}$ such that $\text{Re}\,\beta_2 > -1$, $a,b\in \mathbb{R}$ and $\varepsilon_1, \varepsilon_2 > 0$.*

1. *Suppose $\varepsilon_1 > \varepsilon_2$. Let*

$$I_1(a,b;\beta_1,\beta_2;\varepsilon_1,\varepsilon_2;s,t) = \int_s^t du(-i(u-a)+2\varepsilon_1)^{\beta_1}(-i(u-b)+2\varepsilon_2)^{\beta_2}.$$

*Then*

$$(4.6) \quad \begin{aligned} &I_1(a,b;\beta_1,\beta_2;\varepsilon_1,\varepsilon_2;s,t) \\ &= F_1(a,b;\beta_1,\beta_2;\varepsilon_1,\varepsilon_2;t) - F_1(a,b;\beta_1,\beta_2;\varepsilon_1,\varepsilon_2;s) \\ &= \Phi_1(a,b;\beta_1,\beta_2;\varepsilon_1,\varepsilon_2;t) - \Phi_1(a,b;\beta_1,\beta_2;\varepsilon_1,\varepsilon_2;s), \end{aligned}$$

*with*

$$(4.7) \quad \begin{aligned} &F_1(a,b;\beta_1,\beta_2;\varepsilon_1,\varepsilon_2;t) \\ &= i\frac{(2\varepsilon_2 - i(t-b))^{\beta_2+1}}{\beta_2+1}(2(\varepsilon_1-\varepsilon_2)-i(b-a))^{\beta_1} \\ &\quad \times {}_2F_1\left(-\beta_1,\beta_2+1;\beta_2+2;-\left(\frac{2\varepsilon_2-i(t-b)}{2(\varepsilon_1-\varepsilon_2)-i(b-a)}\right)\right) \end{aligned}$$

*and*

$$(4.8) \quad \begin{aligned} &\Phi_1(a,b;\beta_1,\beta_2;\varepsilon_1,\varepsilon_2;t) \\ &= i\frac{(2\varepsilon_2-i(t-b))^{\beta_1+\beta_2+1}}{\beta_1+\beta_2+1} \\ &\quad \times {}_2F_1\Big(-\beta_1,-(\beta_1+\beta_2+1);-\beta_1-\beta_2; \\ &\qquad\qquad -\left(\frac{2(\varepsilon_1-\varepsilon_2)-i(b-a)}{2\varepsilon_2-i(t-b)}\right)\Big). \end{aligned}$$

2. *(No hypothesis is needed here.) Let*

$$I_2(a,b;\beta_1,\beta_2;\varepsilon_1,\varepsilon_2;s,t) = \int_s^t du(i(u-a)+2\varepsilon_1)^{\beta_1}(-i(u-b)+2\varepsilon_2)^{\beta_2}$$

(notice the sign change with respect to $I_1$). Then

$$
\begin{aligned}
I_2&(a,b;\beta_1,\beta_2;\varepsilon_1,\varepsilon_2;s,t) \\
&= F_2(a,b;\beta_1,\beta_2;\varepsilon_1,\varepsilon_2;t) - F_2(a,b;\beta_1,\beta_2;\varepsilon_1,\varepsilon_2;s),
\end{aligned}
\tag{4.9}
$$

with

$$
\begin{aligned}
F_2&(a,b;\beta_1,\beta_2;\varepsilon_1,\varepsilon_2;t) \\
&= i\frac{(2\varepsilon_2 - i(t-b))^{\beta_2+1}}{\beta_2+1}(2(\varepsilon_1+\varepsilon_2)+i(b-a))^{\beta_1} \\
&\quad \times {}_2F_1\left(-\beta_1,\beta_2+1;\beta_2+2;\left(\frac{2\varepsilon_2 - i(t-b)}{2(\varepsilon_1+\varepsilon_2)+i(b-a)}\right)\right).
\end{aligned}
\tag{4.10}
$$

If, furthermore, $a=b=0$ and $s,t>0$, then

$$
\begin{aligned}
I_2&(0,0;\beta_1,\beta_2;\varepsilon_1,\varepsilon_2;s,t) \\
&= \Phi_2(\beta_1,\beta_2;\varepsilon_1,\varepsilon_2;t) - \Phi_2(\beta_1,\beta_2;\varepsilon_1,\varepsilon_2;s),
\end{aligned}
\tag{4.11}
$$

$$
\begin{aligned}
\Phi_2&(\beta_1,\beta_2;\varepsilon_1,\varepsilon_2;t) \\
&= ie^{i\pi\beta_1}\frac{(2\varepsilon_2-it)^{\beta_1+\beta_2+1}}{\beta_1+\beta_2+1} \\
&\quad \times {}_2F_1\left(-\beta_1,-(\beta_1+\beta_2+1);-\beta_1-\beta_2;\left(\frac{2(\varepsilon_1+\varepsilon_2)}{2\varepsilon_2-it}\right)\right).
\end{aligned}
\tag{4.12}
$$

REMARK. Note the phase factor $e^{i\pi\beta_1}$ in (4.12) is absent from (4.8). It is responsible for the divergence of Lévy's area for $\alpha \leq 1/4$, as will appear clearly during the course of the proof of Theorem 4.4 below.

PROOF OF LEMMA 4.1.

1. Set $v = u - b$. Then, using a deformation of contour,

$$I_1 = F_1(t) - F_1(s)$$

with

$$
F_1(t) = \int_{-2i\varepsilon_2}^{t-b} dv(-iv+2\varepsilon_2)^{\beta_2}(-i(v+b-a)+2\varepsilon_1)^{\beta_1}.
\tag{4.13}
$$

Note that the contour of integration has been implicitly chosen to be a line, though any contour going from $-2i\varepsilon_2$ to $t-b$ while staying above the horizontal line $\operatorname{Im} z = -2\varepsilon_2$ would do; this condition is necessary to avoid problems of multi-valuedness. Note also that $F_1(t)$ may not be defined if $\beta_2 \leq -1$ or if $\varepsilon_1 < \varepsilon_2$ [the integral may diverge, and the function $v \to (i(v+b-a)+2\varepsilon_1)^{\beta_1}$ may not be well defined due to problems of multi-valuedness]. With all this made clear, one may now set $x := \frac{v+2i\varepsilon_2}{t-b+2i\varepsilon_2}$ and obtain

$$F_1(t) = (t - b + 2i\varepsilon_2)$$

(4.14)
$$\times \int_0^1 dx (2\varepsilon_2 - i(t-b))^{\beta_2} x^{\beta_2} (2(\varepsilon_1 - \varepsilon_2) - i(b-a))^{\beta_1}$$

$$\times \left(1 + \left(\frac{2\varepsilon_2 - i(t-b)}{2(\varepsilon_1 - \varepsilon_2) - i(b-a)}\right) x\right)^{\beta_1}$$

$$= i \frac{(2\varepsilon_2 - i(t-b))^{\beta_2+1}}{\beta_2 + 1} (2(\varepsilon_1 - \varepsilon_2) - i(b-a))^{\beta_1}$$

(4.15)
$$\times {}_2F_1\left(-\beta_1, \beta_2 + 1; \beta_2 + 2; -\left(\frac{2\varepsilon_2 - i(t-b)}{2(\varepsilon_1 - \varepsilon_2) - i(b-a)}\right)\right)$$

by formula (4.2).

Now, by the connection formula (4.3), one gets

$$F_1(t) = i \frac{(2\varepsilon_2 - i(t-b))^{\beta_2+1}}{\beta_2 + 1} (2(\varepsilon_1 - \varepsilon_2) - i(b-a))^{\beta_1}$$

$$\times \left\{ \frac{\beta_2 + 1}{\beta_1 + \beta_2 + 1} \left(\frac{2\varepsilon_2 - i(t-b)}{2(\varepsilon_1 - \varepsilon_2) - i(b-a)}\right)^{\beta_1} \right.$$

$$\times {}_2F_1\left(-\beta_1, -(1 + \beta_1 + \beta_2); -\beta_1 - \beta_2; \right.$$

(4.16)
$$\left. -\left(\frac{2(\varepsilon_1 - \varepsilon_2) - i(b-a)}{2\varepsilon_2 - i(t-b)}\right)\right)$$

$$+ \frac{\Gamma(\beta_2 + 2)\Gamma(-(\beta_1 + \beta_2 + 1))}{\Gamma(-\beta_1)}$$

$$\times \left(\frac{2\varepsilon_2 - i(t-b)}{2(\varepsilon_1 - \varepsilon_2) - i(b-a)}\right)^{-1-\beta_2}$$

$$\left. \times {}_2F_1\left(\beta_2 + 1, 0; 2 + \beta_1 + \beta_2; -\left(\frac{2(\varepsilon_1 - \varepsilon_2) - i(b-a)}{2\varepsilon_2 - i(t-b)}\right)\right) \right\}.$$

The second term of the sum simplifies to a constant

$$C_1(a, b; \beta_1, \beta_2; \varepsilon_1, \varepsilon_2)$$
$$= \frac{\Gamma(\beta_2 + 1)\Gamma(-(\beta_1 + \beta_2 + 1))}{\Gamma(-\beta_1)} (2(\varepsilon_1 - \varepsilon_2) - i(b-a))^{\beta_1+\beta_2+1}.$$

By setting $\Phi_1(a, b; \beta_1, \beta_2; \varepsilon_1, \varepsilon_2; s, t) = F_1(a, b; \beta_1, \beta_2; \varepsilon_1, \varepsilon_2; s, t) - C_1(a, b; \beta_1, \beta_2; \varepsilon_1, \varepsilon_2)$, one gets the second formula (4.6).

2. One may use the same deformation of contour, this time without any restriction on the sign of $\varepsilon_1 - \varepsilon_2$, since $\operatorname{Re}(i(v + b - a) + 2\varepsilon_1) > 0$ in any case. Hence, $I_2(s,t) = F_2(t) - F_2(s)$, with

$$F_2(t) = (t - b + 2i\varepsilon_2) \int_0^1 dx (2\varepsilon_2 - i(t - b))^{\beta_2}$$

(4.17)
$$\times x^{\beta_2}(2(\varepsilon_1 + \varepsilon_2) + i(b-a))^{\beta_1}$$

$$\times \left(1 - \left[\frac{2\varepsilon_2 - i(t-b)}{2(\varepsilon_1 + \varepsilon_2) + i(b-a)}\right]x\right)^{\beta_1}$$

and, by (4.2),

$$F_2(a, b; \beta_1, \beta_2; \varepsilon_1, \varepsilon_2; t)$$

(4.18)
$$= i\frac{(2\varepsilon_2 - i(t-b))^{\beta_2+1}}{\beta_2 + 1}(2(\varepsilon_1 + \varepsilon_2) + i(b-a))^{\beta_1}$$

$$\times {}_2F_1\left(-\beta_1, \beta_2 + 1; \beta_2 + 2; \left(\frac{2\varepsilon_2 - i(t-b)}{2(\varepsilon_1 + \varepsilon_2) + i(b-a)}\right)\right).$$

If one wishes to apply formula (4.3) now as above [see (4.16)], one needs to be very careful this time. Calling

$$z(t) = \left(\frac{2\varepsilon_2 - i(t-b)}{2(\varepsilon_1 + \varepsilon_2) + i(b-a)}\right),$$

the argument of the hypergeometric function, it is quite possible that $z(t) \in \mathbb{R}_+$, in which case the connection formula is not defined. However, if $a = b = 0$ and $s, t > 0$, then $z(t), z(s) \in \Pi^-$ and

$$(-z(t))^{\beta_1} = (z(t))^{\beta_1} \cdot e^{i\pi\beta_1}, \qquad (-z(s))^{\beta_1} = (z(s))^{\beta_1} \cdot e^{i\pi\beta_1}.$$

A computation along the same lines as in the case of $I_1$ yields

$$F_2(0, 0; \beta_1, \beta_2; \varepsilon_1, \varepsilon_2; t)$$

(4.19)
$$= C_2(\beta_1, \beta_2; \varepsilon_1, \varepsilon_2) + ie^{i\pi\beta_1}\frac{(2\varepsilon_2 - it)^{\beta_1+\beta_2+1}}{\beta_1 + \beta_2 + 1}$$

$$\times {}_2F_1\left(-\beta_1, -(\beta_1 + \beta_2 + 1); -\beta_1 - \beta_2; \left(\frac{2(\varepsilon_1 + \varepsilon_2)}{2\varepsilon_2 - it}\right)\right),$$

where $C_2(\beta_1, \beta_2; \varepsilon_1, \varepsilon_2)$ does not depend on $t$. □

DEFINITION 4.2 (*Second moment of Lévy's area*). Let

(4.20) $$\mathcal{V}(\varepsilon_1, \varepsilon_2)_t = \mathbb{E}\left[\left(\int_0^t \Gamma'(\varepsilon_1)_{x_1}^{(1)} dx_1 \int_0^{x_1} \Gamma'(\varepsilon_2)_{x_2}^{(2)} dx_2\right)^2\right].$$

If $\varepsilon_1 = \varepsilon_2 =: \varepsilon$, then $\mathcal{V}(\varepsilon_1, \varepsilon_2)_t = \mathbb{E}[(\mathrm{Area}_{0,t}(\mathbf{\Gamma}(\varepsilon)))^2]$.

The more general quantity considered in Definition 4.2 will be needed below for the proof of the convergence of $\mathbf{\Gamma}(\varepsilon)$ when $\varepsilon \to 0$ in the sense of the rough paths.

DEFINITION 4.3. Let $\varepsilon_1, \varepsilon_2 > 0$. We say that $\varepsilon_1, \varepsilon_2 \to 0$ *jointly* if $\varepsilon_1, \varepsilon_2 \to 0$ and $\varepsilon_1/\varepsilon_2 \to 1$.

THEOREM 4.4. *Suppose $\varepsilon_1, \varepsilon_2 \longrightarrow 0$ jointly in the sense of Definition 4.3, and $\alpha > \frac{1}{4}$. Then*

$$(4.21) \qquad \mathcal{V}(\varepsilon_1, \varepsilon_2)_t \longrightarrow C_\alpha t^{4\alpha},$$

*where*

$$(4.22) \quad C_\alpha = \frac{\alpha(2\alpha - 1)}{2} \left[ \frac{2\Gamma(2\alpha - 1)\Gamma(2\alpha + 1)}{\Gamma(4\alpha + 1)} + \frac{1}{(2\alpha - 1)(4\alpha - 1)} \right].$$

PROOF. One has

$$\begin{aligned}
\mathcal{V}(\varepsilon_1, \varepsilon_2)_t &= 2\operatorname{Re}\mathbb{E}\Bigg[\Bigg(\int_0^t \Gamma'(\varepsilon_1)^{+,(1)}_{x_1}\, dx_1 \int_0^{x_1} \Gamma'(\varepsilon_2)^{+,(2)}_{x_2}\, dx_2\Bigg) \\
&\qquad\qquad\times \Bigg(\int_0^t \Gamma'(\varepsilon_1)^{-,(1)}_{y_1}\, dy_1 \int_0^{y_1} \Gamma'(\varepsilon_2)^{-,(2)}_{y_2}\Bigg)\, dy_2\Bigg] \\
(4.23) \qquad &\quad + 2\operatorname{Re}\mathbb{E}\Bigg[\Bigg(\int_0^t \Gamma'(\varepsilon_1)^{+,(1)}_{x_1}\, dx_1 \int_0^{x_1} \Gamma'(\varepsilon_2)^{-,(2)}_{x_2}\, dx_2\Bigg) \\
&\qquad\qquad\times \Bigg(\int_0^t \Gamma'(\varepsilon_1)^{-,(1)}_{y_1}\, dy_1 \int_0^{y_1} \Gamma'(\varepsilon_2)^{+,(2)}_{y_2}\Bigg)\, dy_2\Bigg] \\
&= \left(\frac{\alpha(1 - 2\alpha)}{2\cos \pi\alpha}\right)^2 \cdot 2\operatorname{Re}(\mathcal{V}_1(\varepsilon_1, \varepsilon_2)_t + \mathcal{V}_2(\varepsilon_1, \varepsilon_2)_t),
\end{aligned}$$

where

$$\begin{aligned}
\mathcal{V}_1 &= \int_0^t dx_1 \int_0^t dy_1 (-i(x_1 - y_1) + 2\varepsilon_1)^{2\alpha - 2} \\
&\quad \times \int_0^{x_1} dx_2 \int_0^{y_1} dy_2 (-i(x_2 - y_2) + 2\varepsilon_2)^{2\alpha - 2} \\
(4.24) \qquad &= \frac{1}{2\alpha(2\alpha - 1)} \int_0^t dx_1 \int_0^t dy_1 (-i(x_1 - y_1) + 2\varepsilon_1)^{2\alpha - 2} \\
&\quad \times ((-i(x_1 - y_1) + 2\varepsilon_2)^{2\alpha} - (-ix_1 + 2\varepsilon_2)^{2\alpha} \\
&\qquad\qquad - (iy_1 + 2\varepsilon_2)^{2\alpha} + (2\varepsilon_2)^{2\alpha})
\end{aligned}$$

and

$$\mathcal{V}_2 = \int_0^t dx_1 \int_0^t dy_1(-i(x_1 - y_1) + 2\varepsilon_1)^{2\alpha-2}$$
$$\times \int_0^{x_1} dx_2 \int_0^{y_1} dy_2(i(x_2 - y_2) + 2\varepsilon_2)^{2\alpha-2}$$

(4.25)
$$= \frac{1}{2\alpha(2\alpha - 1)} \int_0^t dx_1 \int_0^t dy_1(-i(x_1 - y_1) + 2\varepsilon_1)^{2\alpha-2}$$
$$\times ((i(x_1 - y_1) + 2\varepsilon_2)^{2\alpha} - (ix_1 + 2\varepsilon_2)^{2\alpha}$$
$$- (-iy_1 + 2\varepsilon_2)^{2\alpha} + (2\varepsilon_2)^{2\alpha}).$$

*Computation of $\mathcal{V}_1$.* From the last expression (4.24), one gets

(4.26)
$$\mathcal{V}_1 = \frac{1}{2\alpha(2\alpha - 1)}[\mathcal{V}_{1,1} - \mathcal{V}_{1,2} - \mathcal{V}_{1,3} + \mathcal{V}_{1,4}],$$

where

$$\mathcal{V}_{1,1} = \int_0^t dx_1 \int_0^t dy_1(-i(x_1 - y_1) + 2\varepsilon_1)^{2\alpha-2}(-i(x_1 - y_1) + 2\varepsilon_2)^{2\alpha}$$

(4.27)
$$= \int_{-t}^t dx(-ix + 2\varepsilon_1)^{2\alpha-2}(-ix + 2\varepsilon_2)^{2\alpha}(t - |x|)$$
$$= 2\operatorname{Re}(tI_1(0, 0; 2\alpha - 2, 2\alpha; \varepsilon_1, \varepsilon_2; 0, t) - \mathcal{V}_{1,1,2}),$$

where

$$\mathcal{V}_{1,1,2} = \int_0^t dx\, x(-ix + 2\varepsilon_1)^{2\alpha-2}(-ix + 2\varepsilon_2)^{2\alpha}$$
$$= \frac{i}{2\alpha - 1}\Big\{(-it + 2\varepsilon_1)^{2\alpha-1}t(-it + 2\varepsilon_2)^{2\alpha}$$
$$- I_1(0, 0; 2\alpha - 1, 2\alpha; \varepsilon_1, \varepsilon_2; 0, t)$$
$$+ 2i\alpha \int_0^t dx((x + i(\varepsilon_1 + \varepsilon_2)) - i(\varepsilon_1 + \varepsilon_2))$$
$$\times ((-ix + 2\varepsilon_1)(-ix + 2\varepsilon_2))^{2\alpha-1}\Big\}$$
$$= \frac{i}{2\alpha - 1}\Big\{(-it + 2\varepsilon_1)^{2\alpha-1}t(-it + 2\varepsilon_2)^{2\alpha}$$
$$- I_1(0, 0; 2\alpha - 1, 2\alpha; \varepsilon_1, \varepsilon_2; 0, t)$$
$$- \frac{i}{2}[((-it + 2\varepsilon_1)(-it + 2\varepsilon_2))^{2\alpha} - (4\varepsilon_1\varepsilon_2)^{2\alpha}]$$
$$+ 2\alpha(\varepsilon_1 + \varepsilon_2)I_1(0, 0; 2\alpha - 1, 2\alpha - 1; \varepsilon_1, \varepsilon_2; 0, t)\Big\}.$$

Suppose $\varepsilon_1, \varepsilon_2 \to 0$ jointly (in the sense of Definition 4.3). Then [by using the expression (4.6)]

$$\operatorname{Re} I_1(0,0; 2\alpha - 2, 2\alpha; \varepsilon_1, \varepsilon_2; 0, t)$$

(4.28)
$$= \frac{1}{4\alpha - 1} \operatorname{Re} i \bigg\{ (2\varepsilon_2 - it)^{4\alpha - 1} {}_2F_1\bigg(2 - 2\alpha, 1 - 4\alpha; 2 - 4\alpha;$$
$$- \bigg(\frac{2(\varepsilon_1 - \varepsilon_2)}{2\varepsilon_2 - it}\bigg)\bigg)\bigg\}$$
$$\longrightarrow \frac{1}{4\alpha - 1} \operatorname{Re} i(-it)^{4\alpha - 1} = -\frac{1}{4\alpha - 1} \cos 2\pi\alpha \, t^{4\alpha - 1}.$$

[Note that $\operatorname{Im} I_1(0,0; 2\alpha - 2, 2\alpha; \varepsilon_1, \varepsilon_2; 0, t)$ diverges in the same limit!]

Next, by (4.6),

(4.29)
$$I_1(0,0; 2\alpha - 1, 2\alpha - 1; \varepsilon_1, \varepsilon_2; 0, t)$$
$$= \frac{i}{4\alpha - 1} \bigg\{ (2\varepsilon_2 - it)^{4\alpha - 1}$$
$$\times {}_2F_1\bigg(1 - 2\alpha, 1 - 4\alpha; 2 - 4\alpha; -\bigg(\frac{2(\varepsilon_1 - \varepsilon_2)}{2\varepsilon_2 - it}\bigg)\bigg)$$
$$- (2\varepsilon_2)^{4\alpha - 1} {}_2F_1\bigg(1 - 2\alpha, 1 - 4\alpha; 2 - 4\alpha;$$
$$- \bigg(\frac{2(\varepsilon_1 - \varepsilon_2)}{2\varepsilon_2}\bigg)\bigg)\bigg\},$$

hence, $(\varepsilon_1 + \varepsilon_2) I_1(0,0; 2\alpha - 1, 2\alpha - 1; \varepsilon_1, \varepsilon_2; 0, t) \to 0$. Note that this expression has no limit if $\varepsilon_1$ and $\varepsilon_2$ go independently to 0, and not jointly (in the sense of Definition 4.3).

Finally,

(4.30)
$$I_1(0,0; 2\alpha - 1, 2\alpha; \varepsilon_1, \varepsilon_2; 0, t)$$
$$= \frac{i(2\varepsilon_2 - it)^{4\alpha}}{4\alpha} {}_2F_1\bigg(1 - 2\alpha, -4\alpha; 1 - 4\alpha; -\bigg(\frac{2(\varepsilon_1 - \varepsilon_2)}{2\varepsilon_2 - it}\bigg)\bigg)$$
$$- \frac{i(2\varepsilon_2)^{4\alpha}}{4\alpha} {}_2F_1\bigg(1 - 2\alpha, -4\alpha; 1 - 4\alpha; -\bigg(\frac{2(\varepsilon_1 - \varepsilon_2)}{2\varepsilon_2}\bigg)\bigg)$$
$$\to i \frac{(-it)^{4\alpha}}{4\alpha} = i e^{-2i\pi\alpha} \frac{t^{4\alpha}}{4\alpha}.$$

Hence,

(4.31)
$$\mathcal{V}_{1,1,2} \to -\frac{1}{4\alpha} e^{-2i\pi\alpha} t^{4\alpha}$$

and

(4.32) $$\mathcal{V}_{1,1} \to -\frac{\cos 2\pi\alpha}{2\alpha(4\alpha-1)}t^{4\alpha}.$$

Let us now turn to the computation of $\mathcal{V}_{1,2}$:

$$\mathcal{V}_{1,2} = \int_0^t dx_1 \int_0^t dy_1 (-i(x_1-y_1)+2\varepsilon_1)^{2\alpha-2}(-ix_1+2\varepsilon_2)^{2\alpha}$$

$$= -\frac{i}{2\alpha-1}\int_0^t dx_1(-ix_1+2\varepsilon_2)^{2\alpha}$$

(4.33)
$$\times [(-i(x_1-t)+2\varepsilon_1)^{2\alpha-1} - (-ix_1+2\varepsilon_1)^{2\alpha-1}]$$

$$= -\frac{i}{2\alpha-1}(I_1(t,0;2\alpha-1,2\alpha;\varepsilon_1,\varepsilon_2;0,t)$$

$$- I_1(0,0;2\alpha-1,2\alpha;\varepsilon_1,\varepsilon_2;0,t)).$$

Using (4.6) this time,

$$I_1(t,0;2\alpha-1,2\alpha;\varepsilon_1,\varepsilon_2;0,t)$$

$$= \frac{i}{2\alpha+1}\bigg\{(2\varepsilon_2-it)^{2\alpha+1}(2(\varepsilon_1-\varepsilon_2)+it)^{2\alpha-1}$$

$$\times {}_2F_1\bigg(1-2\alpha,1+2\alpha;2+2\alpha;-\bigg(\frac{2\varepsilon_2-it}{2(\varepsilon_1-\varepsilon_2)+it}\bigg)\bigg)$$

(4.34)
$$- (2\varepsilon_2)^{2\alpha+1}(2(\varepsilon_1-\varepsilon_2)+it)^{2\alpha-1}$$

$$\times {}_2F_1\bigg(1-2\alpha,1+2\alpha;2+2\alpha;-\bigg(\frac{2\varepsilon_2}{2(\varepsilon_1-\varepsilon_2)+it}\bigg)\bigg)\bigg\}$$

$$\longrightarrow \frac{i}{2\alpha+1}(-it)^{2\alpha+1}(it)^{2\alpha-1}{}_2F_1(1-2\alpha,1+2\alpha;2+2\alpha;1)$$

$$= -\frac{i}{2\alpha}\frac{\Gamma(2\alpha+1)^2}{\Gamma(4\alpha+1)}t^{4\alpha}$$

by (4.5).

The function $I_1(0,0;2\alpha-1,2\alpha;\varepsilon_1,\varepsilon_2;0,t)$ has already been studied; see (4.30). Hence,

(4.35) $$\mathcal{V}_{1,2} \to -\frac{1}{2\alpha-1}\bigg(\frac{\Gamma(2\alpha)\Gamma(2\alpha+1)}{\Gamma(4\alpha+1)} + \frac{e^{-2i\pi\alpha}}{4\alpha}\bigg)t^{4\alpha}.$$

Finally, $\mathcal{V}_{1,3} = \overline{\mathcal{V}_{1,2}}$ and, clearly,

$$\mathcal{V}_{1,4} = \bigg(\frac{2\cos\pi\alpha}{\alpha(1-2\alpha)}\mathbb{E}[\Gamma(\varepsilon_1)_t^+\Gamma(\varepsilon_1)_t^-]\bigg)(2\varepsilon_2)^{2\alpha} \stackrel{\varepsilon_1,\varepsilon_2\to 0}{\sim} \frac{2\cos\pi\alpha}{\alpha(1-2\alpha)}\frac{t^{2\alpha}}{2}(2\varepsilon_2)^{2\alpha}$$

$$\to 0.$$

Hence,

$$\mathcal{V}_1(\varepsilon_1, \varepsilon_2)_t \longrightarrow \frac{1}{2\alpha(2\alpha - 1)}$$
$$\times \left[ \frac{2\Gamma(2\alpha - 1)\Gamma(2\alpha + 1)}{\Gamma(4\alpha + 1)} + \frac{\cos 2\pi\alpha}{(2\alpha - 1)(4\alpha - 1)} \right] t^{4\alpha}. \quad (4.36)$$

One may check by using the asymptotic expansion of the Gamma function near zero,

$$\Gamma(\varepsilon) \sim_{\varepsilon \to 0} \frac{1}{\varepsilon} + \gamma + O(\varepsilon)$$

(where $\gamma$ is Euler's constant), that this expression is regular when $\alpha \to \frac{1}{4}$ or $\alpha \to \frac{1}{2}$.

*Computation of $\mathcal{V}_2$.* The computations are very similar to the previous ones, but the final result is surprisingly different, owing to the presence of unfortunate (but unevitable) phase factors coming from point 2 of Lemma 4.1. The reader may easily check step by step that (keeping the same notation)

$$\mathcal{V}_2 = \frac{1}{2\alpha(2\alpha - 1)} (\mathcal{V}_{2,1} - \mathcal{V}_{2,2} - \mathcal{V}_{2,3} + \mathcal{V}_{2,4});$$

$$\mathcal{V}_{2,1} = 2\operatorname{Re}(t\overline{I_2(0, 0; 2\alpha - 2, 2\alpha; \varepsilon_1, \varepsilon_2; 0, t)} - \mathcal{V}_{2,1,2}), \quad (4.37)$$

where

$$\mathcal{V}_{2,1,2} = \int_0^t dx\, x(-ix + 2\varepsilon_1)^{2\alpha - 2}(ix + 2\varepsilon_2)^{2\alpha}$$

$$= \frac{i}{2\alpha - 1} \Big\{ (-it + 2\varepsilon_1)^{2\alpha - 1} t(it + 2\varepsilon_2)^{2\alpha}$$

$$- \overline{I_2(0, 0; 2\alpha - 1, 2\alpha; \varepsilon_1, \varepsilon_2; 0, t)}$$

$$- 2i\alpha \int_0^t dx((x + i(\varepsilon_1 - \varepsilon_2)) - i(\varepsilon_1 - \varepsilon_2))$$

$$\times ((-ix + 2\varepsilon_1)(ix + 2\varepsilon_2))^{2\alpha - 1} \Big\}$$

$$= \frac{i}{2\alpha - 1} \Big\{ (-it + 2\varepsilon_1)^{2\alpha - 1} t(it + 2\varepsilon_2)^{2\alpha}$$

$$- \overline{I_2(0, 0; 2\alpha - 1, 2\alpha; \varepsilon_1, \varepsilon_2; 0, t)}$$

$$- \frac{i}{2} [((-it + 2\varepsilon_1)(it + 2\varepsilon_2))^{2\alpha} - (4\varepsilon_1 \varepsilon_2)^{2\alpha}]$$

$$- 2\alpha(\varepsilon_1 - \varepsilon_2)\overline{I_2(0, 0; 2\alpha - 1, 2\alpha - 1; \varepsilon_1, \varepsilon_2; 0, t)} \Big\}.$$

Suppose $\varepsilon_1, \varepsilon_2 \to 0$ jointly. Then [by using expression (4.11)]

$$I_2(0,0; 2\alpha - 2, 2\alpha; \varepsilon_1, \varepsilon_2; 0, t)$$

$$= \frac{1}{4\alpha - 1} i e^{2i\pi\alpha} (2\varepsilon_2 - it)^{4\alpha - 1}$$

(4.38)
$$\times {}_2F_1\left(2 - 2\alpha, 1 - 4\alpha; 2 - 4\alpha; \left(\frac{2(\varepsilon_1 + \varepsilon_2)}{2\varepsilon_2 - it}\right)\right)$$

$$- \frac{1}{4\alpha - 1} i e^{2i\pi\alpha} (2\varepsilon_2)^{4\alpha - 1}$$

$$\times \lim_{s \overset{>}{\to} 0} {}_2F_1\left(2 - 2\alpha, 1 - 4\alpha; 2 - 4\alpha; \left(\frac{2(\varepsilon_1 + \varepsilon_2)}{2\varepsilon_2 - is}\right)\right).$$

The first term in the right-hand side goes to

$$-\frac{1}{4\alpha - 1} t^{4\alpha - 1}$$

as $\varepsilon_1, \varepsilon_2 \to 0$, but the second term is equivalent as $\varepsilon_1, \varepsilon_2 \to 0$ while $\varepsilon_1/\varepsilon_2 \to 1$ to

(4.39) $\quad -\frac{1}{4\alpha - 1} i e^{2i\pi\alpha} (2\varepsilon_2)^{4\alpha - 1} \lim_{s \overset{>}{\to} 0} {}_2F_1(2 - 2\alpha, 1 - 4\alpha; 2 - 4\alpha; 2 + is),$

whose real part does not miraculously vanish this time, because of the phase factor $e^{2i\pi\alpha}$. However, *if* $\alpha > \frac{1}{4}$, then this goes to zero.

The other terms prove to be regular for any $\alpha > 0$: namely,

(4.40) $\quad (\varepsilon_1 - \varepsilon_2) I_2(0, 0; 2\alpha - 1, 2\alpha - 1; \varepsilon_1, \varepsilon_2; 0, t) \longrightarrow 0;$

(4.41)
$$I_2(0, 0; 2\alpha - 1, 2\alpha; \varepsilon_1, \varepsilon_2; 0, t) \to i e^{i\pi(2\alpha - 1)} \frac{(-it)^{4\alpha}}{4\alpha}$$

$$= -i \frac{t^{4\alpha}}{4\alpha}$$

by (4.12).

Hence,

(4.42) $\quad \mathcal{V}_{2,1,2} \to -\frac{1}{4\alpha} t^{4\alpha}$

and, supposing $\alpha > \frac{1}{4}$,

(4.43) $\quad \mathcal{V}_{2,1} \to -\frac{1}{2\alpha(4\alpha - 1)} t^{4\alpha}.$

Let us now turn to the computation of $\mathcal{V}_{2,2}$:

(4.44)
$$\mathcal{V}_{2,2} = -\frac{i}{2\alpha - 1} (\overline{I_2(t, 0; 2\alpha - 1, 2\alpha; \varepsilon_1, \varepsilon_2; 0, t)}$$

$$- \overline{I_2(0, 0; 2\alpha - 1, 2\alpha; \varepsilon_1, \varepsilon_2; 0, t)});$$

using (4.9) this time,

$$I_2(t, 0; 2\alpha - 1, 2\alpha; \varepsilon_1, \varepsilon_2; 0, t)$$

$$= \frac{i}{2\alpha + 1}\Big\{ (2\varepsilon_2 - it)^{2\alpha+1}(2(\varepsilon_1 + \varepsilon_2) - it)^{2\alpha-1}$$

$$\times {}_2F_1\Big(1 - 2\alpha, 1 + 2\alpha; 2 + 2\alpha; \Big(\frac{2\varepsilon_2 - it}{2(\varepsilon_1 + \varepsilon_2) - it}\Big)\Big)$$

(4.45)
$$- (2\varepsilon_2)^{2\alpha+1}(2(\varepsilon_1 + \varepsilon_2) - it)^{2\alpha-1}$$

$$\times {}_2F_1\Big(1 - 2\alpha, 1 + 2\alpha; 2 + 2\alpha; \Big(\frac{2\varepsilon_2}{2(\varepsilon_1 + \varepsilon_2) - it}\Big)\Big)\Big\}$$

$$\longrightarrow \frac{i}{2\alpha + 1}(-it)^{4\alpha} {}_2F_1(1 - 2\alpha, 1 + 2\alpha; 2 + 2\alpha; 1)$$

$$= \frac{i}{2\alpha}\frac{\Gamma(2\alpha + 1)^2}{\Gamma(4\alpha + 1)} e^{-2i\pi\alpha} t^{4\alpha}$$

by (4.5).

Hence,

(4.46) $$\mathcal{V}_{2,2} \longrightarrow -\frac{1}{2\alpha - 1}\Big(\frac{\Gamma(2\alpha)\Gamma(2\alpha + 1)}{\Gamma(4\alpha + 1)} e^{2i\pi\alpha} + \frac{1}{4\alpha}\Big) t^{4\alpha}$$

and, if $\alpha > \frac{1}{4}$,

(4.47)
$$\mathcal{V}_2 \longrightarrow \frac{1}{2\alpha(2\alpha - 1)}$$

$$\times \Big[2\cos 2\pi\alpha \frac{\Gamma(2\alpha - 1)\Gamma(2\alpha + 1)}{\Gamma(4\alpha + 1)} + \frac{1}{(2\alpha - 1)(4\alpha - 1)}\Big] t^{4\alpha}.$$

Hence, if $\alpha > \frac{1}{4}$, then (after some simplifications) one gets (4.22). □

REMARKS. Let us compute the constant $C_\alpha$ in some particular cases.
When $\alpha \to 1$, $C_\alpha \to \frac{1}{4}$.
When $\alpha \to \frac{1}{2}$, $C_\alpha \to \frac{1}{2}$ [note that the apparent singularity in formula (4.22) vanishes].
When $\alpha \overset{>}{\to} \frac{1}{4}$, $C_\alpha \sim \frac{1/8}{4\alpha - 1} \to +\infty$.

Let us conclude this paragraph by mentioning that the piecewise linear approximation studied in [5] leads to the same value for $C_\alpha$ (in other words, for the Lévy area of the limiting process) as proved in a recent work by Baudoin and Coutin (see [3]), although a few integration by parts are required to see that. Hence, a natural question arises: do both approximations define in the limit the same stochastic integration theory? We have no answer for the moment, although we believe this should be true.

Also, let $\mathcal{A}_t := \int_0^t B_s^{(2)}\, dB_s^{(1)} - B(1)_s\, dB_s^{(2)}\, ds$ be the usual antisymmetric Lévy area for FBM. Then an easy integration by parts yields $\mathbb{E}[\mathcal{A}_t^2] = (4C_\alpha - 1)t^{4\alpha}$.

4.2. *Convergence of Lévy's area for $\alpha > \frac{1}{4}$.* We shall prove in this subsection the $d_2$-convergence of the $q$-variation distance $d_1$, $d_2$ of the Lévy area for $\alpha > \frac{1}{4}$ and $q > \frac{1}{\alpha}$.

We shall be using the following technical lemmas (see [5] and [13]):

PROPOSITION 4.5. *Let $0 < s < t$.*

1. *Let $\mathbf{w}, \mathbf{v}$ be two $V$-valued paths with $q$-bounded variation. Let $t_k^n = s + (t-s)\frac{k}{2^n}$, $k = 0, \ldots, 2^n$, be the nth dyadic subdivision of $[s, t]$. Then for every $\kappa > q/2 - 1$, there exists a constant $C$ depending only on $q$ and $\kappa$ such that*

$$d_1(\mathbf{w}, \mathbf{v})_{s,t} \leq C \left[ \sum_{n \geq 1} n^\kappa \sum_{k=1}^{2^n} |\mathbf{w}^1_{t^n_{k-1}, t^n_k} - \mathbf{v}^1_{t^n_{k-1}, t^n_k}|^q \right]^{1/q}.$$

2. *Let $\mathbf{w}, \mathbf{v}$ be two $V$-valued paths with $q$-bounded variation and bounded $d_2$-norm. Then for every $\kappa > q/2 - 1$, there exists a constant $C$ depending only on $q$ and $\kappa$ such that*

(4.48)
$$\begin{aligned}
(d_2(\mathbf{w}, \mathbf{v})_{s,t})^{q/2} &\leq C \sum_{n \geq 1} n^\kappa \sum_{k=1}^{2^n} |\mathbf{w}^2_{t^n_{k-1}, t^n_k} - \mathbf{v}^2_{t^n_{k-1}, t^n_k}|^{q/2} \\
&\quad + C \left( \sum_{n \geq 1} n^\kappa \sum_{k=1}^{2^n} |\mathbf{w}^1_{t^n_{k-1}, t^n_k} - \mathbf{v}^1_{t^n_{k-1}, t^n_k}|^q \right)^{1/2} \\
&\quad \times \left( \sum_{n \geq 1} n^\kappa \sum_{k=1}^{2^n} |\mathbf{w}^1_{t^n_{k-1}, t^n_k}|^q + |\mathbf{v}^1_{t^n_{k-1}, t^n_k}|^q \right)^{1/2}.
\end{aligned}$$

LEMMA 4.6. *Let*

(4.49) $$\delta(\kappa; \varepsilon; \alpha_1, \alpha_2) = \left[ \sum_{n=0}^{E(|\log_2 \varepsilon|)} n^\kappa 2^n (\varepsilon^{\alpha_1} 2^{-n\alpha_2})^{q/2} \right]^{2/q}$$

*($E$ = entire part) for $\kappa, \varepsilon, \alpha_1 > 0$, with $\alpha_2 := 2\alpha - \alpha_1$. Then*

(4.50) $$\sum_{j \geq 1} \delta(\kappa; j^{-\beta}; \alpha_1, \alpha_2) < \infty,$$

*if $2\beta(\alpha - \frac{1}{q}) > 1$ and $\beta\alpha_1 > 1$.*

Proof. One has

$$\delta(\kappa; \varepsilon; \alpha_1, \alpha_2) = \left[ \sum_{n=0}^{E(|\log_2 \varepsilon|)} n^\kappa 2^{n(1-\alpha_2 q/2)} \right]^{2/q} \cdot \varepsilon^{\alpha_1}.$$

Let $\gamma = (1 - \alpha_2 q/2) \log 2$; the precise estimates depend on the sign of $\gamma$.

Suppose first that $\gamma \leq 0$: then

$$\delta(\kappa; \varepsilon; \alpha_1, \alpha_2) \leq C(1 + |\log_2 \varepsilon|)^{2(\kappa+1)/q} \cdot \varepsilon^{\alpha_1},$$

hence, $\sum_{j \geq 1} \delta(\kappa; j^{-\beta}; \alpha_1, \alpha_2)$ converges if $\beta \alpha_1 > 1$.

On the other hand, if $\gamma > 0$, then

$$\sum_{n=0}^{E(|\log_2 \varepsilon|)} n^\kappa 2^{n(1-\alpha_2 q/2)} \leq \int_0^{|\log_2 \varepsilon|+1} t^\kappa e^{\gamma t}\, dt,$$

hence, there exists $C > 0$ (depending on $q$ and $\kappa$) such that

(4.51)
$$\delta(\kappa; \varepsilon; \alpha_1, \alpha_2) \leq C[(1 + |\log_2 \varepsilon|)^\kappa \varepsilon^{-1+\alpha_2 q/2}]^{2/q} \varepsilon^{\alpha_1}$$
$$= C(1 + |\log_2 \varepsilon|)^{2\kappa/q} \varepsilon^{2(\alpha - 1/q)}$$

and

$$\sum_{j \geq 1} \delta(\kappa; j^{-\beta}; \alpha_1, \alpha_2) \leq C \sum_{j \geq 1} (1 + \beta \log j)^{2\kappa/q} j^{-2\beta(\alpha - 1/q)},$$

which converges if $2\beta(\alpha - 1/q) > 1$. □

LEMMA 4.7. *Let $t_k^n$, $k = 0, \ldots, 2^n$, be the nth dyadic partition of $[0,1]$. Set, for $\varepsilon > \eta > 0$ and $\kappa > 0$,*

(4.52)
$$\delta'(\kappa; d; \varepsilon, \eta)$$
$$= \left[ \sum_{n \geq E(|\log_2 \eta|)} n^\kappa \sum_{k=1}^n \|\mathbf{\Gamma}(\varepsilon)_{t_{k-1}^n, t_k^n}^d - \mathbf{\Gamma}(\eta)_{t_{k-1}^n, t_k^n}^d\|^{q/d} \right]^{d/q},$$

*where $\mathbf{\Gamma}(\varepsilon)^d \in V^{\otimes d}$ is the tensor d-iterated integral; see (0.4).*

*Then with probability one,*

$$\sum_{j \geq 1} \delta'(\kappa; d; j^{-\beta}, (j+1)^{-\beta}) < \infty,$$

*if $d\beta(\alpha - 1/q) > 1$.*

PROOF. Let, for any choice of indices $i_1, \ldots, i_d$,

$$
\begin{aligned}
I_d(\varepsilon, \eta; t) \\
= \mathbb{E}\bigg[\bigg(&\int_{t_{k-1}^n}^{t_k^n} \Gamma'(\varepsilon)_{x_1}^{(i_1)} dx_1 \int_{t_{k-1}^n}^{x_1} \Gamma'(\varepsilon)_{x_2}^{(i_2)} dx_2 \cdots \int_{t_{k-1}^n}^{x_{d-1}} \Gamma'(\varepsilon)_{x_d}^{(i_d)} dx_d \\
&- \int_{t_{k-1}^n}^{t_k^n} \Gamma'(\eta)_{x_1}^{(i_1)} dx_1 \int_{t_{k-1}^n}^{x_1} \Gamma'(\eta)_{x_2}^{(i_2)} dx_2 \\
&\qquad \cdots \int_{t_{k-1}^n}^{x_{d-1}} dx_d \Gamma'(\eta)_{x_d}^{(i_d)} dx_d\bigg)^2\bigg].
\end{aligned}
\tag{4.53}
$$

Then one may estimate $I_d$ very crudely by

$$
\begin{aligned}
I_d(\varepsilon, \eta; t) &\leq C\bigg[\int_0^{2^{-n}} dx \int_0^{2^{-n}} dy |-i(x-y) + 2\eta|^{2\alpha-2}\bigg]^d \\
&\leq C(2^{-2n} \eta^{2\alpha-2})^d.
\end{aligned}
\tag{4.54}
$$

So

$$
\begin{aligned}
\mathbb{E}[\delta'(\kappa; \varepsilon, \eta)] &\leq C'\bigg(\sum_{n > E(|\log_2 \eta|)} n^\kappa 2^n (2^{-2n}\eta^{2\alpha-2})^{q/2}\bigg)^{d/q} \\
&= C'\bigg(\sum_{n > E(|\log_2 \eta|)} n^\kappa 2^{n(1-q)}\bigg)^{d/q} \eta^{d(\alpha-1)} \\
&\leq C'\bigg(\int_{|\log_2 \eta|}^\infty t^\kappa e^{-t(q-1)\log 2} dt\bigg)^{d/q} \eta^{d(\alpha-1)} \\
&\leq C''(1 + |\log_2 \eta|)^{\kappa d/q} \eta^{d(\alpha-1/q)},
\end{aligned}
\tag{4.55}
$$

hence the result. $\square$

COROLLARY 4.8. *Suppose there exists a constant $C > 0$ such that for every $s, t \in \mathbb{R}$, $\varepsilon > \eta > 0$ with $|s - t| \geq C\varepsilon$,*

$$
\mathbb{E}[\|\mathbf{\Gamma}(\varepsilon)_{s,t}^2 - \mathbf{\Gamma}(\eta)_{s,t}^2\|^2] \leq \sum_{i=1}^I C_i \varepsilon^{2\alpha_{i,1}} |t-s|^{2\alpha_{i,2}}
\tag{4.56}
$$

*for some constants $C_i > 0$ and couples of exponents $(\alpha_{i_1}, \alpha_{i_2})$ such that $\alpha_{i_1} + \alpha_{i_2} = 2\alpha$ and $\alpha_{i_1} > 0$.*

*Then the series $\sum_{j \geq 1}(\mathbf{\Gamma}(j^{-\beta})_{s,t}^2 - \mathbf{\Gamma}((j+1)^{-\beta})_{s,t}^2), j \geq 1)$ converges a.s. in the $d_2$-norm if $2\beta(\alpha - 1/q) > 1$ and $\beta\alpha_{i,1} > 1$, $i = 1, \ldots, I$. Hence, the $\varepsilon$-regularized Lévy area $\mathbf{\Gamma}(\varepsilon)_{s,t}^2$ converges a.s. in the $d_2$-norm when $\varepsilon \to 0$.*

PROOF. We shall actually prove that $\mathbb{E}[\sum_{j\geq 1} d_2(\mathbf{\Gamma}(j^{-\beta})_{0,t}, \mathbf{\Gamma}((j+1)^{-\beta})_{0,t})]$ converges. Namely, by Proposition 4.5,

$$
\begin{aligned}
(4.57) \quad & d_2^{q/2}(\mathbf{\Gamma}(j^{-\beta}), \mathbf{\Gamma}((j+1)^{-\beta}))_{0,t} \\
& \leq \left( \sum_{n=0}^{E(\beta \log_2 j)} + \sum_{n > E(\beta \log_2 j)} \right) \\
& \quad \times n^\kappa \sum_{k=1}^{2^n} \|\mathbf{\Gamma}(j^{-\beta})^2_{t^n_{k-1}, t^n_k} - \mathbf{\Gamma}((j+1)^{-\beta})^2_{t^n_{k-1}, t^n_k}\|^{q/2}.
\end{aligned}
$$

It is a well-known consequence of the hypercontractivity property of the Ornstein–Uhlenbeck process (see [16]) that $L^p$-norms ($1 < p < \infty$) are all equivalent on any Wiener chaos $\mathcal{H}_n$ (in this case $n = 2$). Hence, the first term in the right-hand side of the last equation is bounded by

$$
\begin{aligned}
(4.58) \quad & C \sum_{i=1}^{I} C_i^{q/4} \sum_{n=0}^{[\beta \log_2 j]} n^\kappa 2^n (2^{-n})^{q\alpha_{i,2}/2} (j^{-\beta})^{q\alpha_{i,1}/2} \\
& \leq C' \sum_{i=1}^{I} (\delta(\kappa; j^{-\beta}; \alpha_{i_1}, \alpha_{i_2}))^{q/2}
\end{aligned}
$$

(see Lemma 4.6).

On the other hand, the second term is bounded by $(\delta'(\kappa; 2; j^{-\beta}, (j+1)^{-\beta}))^{q/2}$ (see Lemma 4.7).

Hence (by Hölder's inequality, provided $q/2 > 1$, which is the case when $\alpha < \frac{1}{2}$),

$$
\begin{aligned}
& \mathbb{E}\left[ \sum_{j \geq 1} d_2(\mathbf{\Gamma}(j^{-\beta})_{0,t}, \mathbf{\Gamma}(j^{-\beta+1})_{0,t}) \right] \\
& \leq C \sum_{i=1}^{I} \left( \sum_{j \geq 1} \delta(\kappa; j^{-\beta}; \alpha_{i_1}, \alpha_{i_2}) \right) + \sum_{j \geq 1} \delta'(\kappa; 2; j^{-\beta}, (j+1)^{-\beta}) < \infty
\end{aligned}
$$

(by Lemmas 4.6 and 4.7) if $2\beta(\alpha - 1/q) > 1$ and $\beta \alpha_{i,1} > 1$, $i = 1, \ldots, I$.

The convergence of the regularized Lévy area follows easily along the same lines by completeness of the space of rough paths $\Omega_q(V)$. □

We shall need below the following definition of the *second variation* of a function with respect to one variable:

DEFINITION 4.9. Let $f : \mathbb{R}^n \to \mathbb{R}$ be any function of $n$ variables. Then the second variation of $f$ with respect to the $k$th variable is by definition the following function of $(n+1)$ variables:

$$\Delta_2(f(x_1,\ldots,x_{k-1},\ldots,x_{k+1},\ldots,x_n))(\varepsilon,\eta)$$

(4.59)
$$:= f(x_1,\ldots,x_{k-1},\varepsilon,x_{k+1},\ldots,x_n)$$
$$- 2f\left(x_1,\ldots,x_{k-1},\frac{\varepsilon+\eta}{2},x_{k+1},\ldots,x_n\right)$$
$$+ f(x_1,\ldots,x_{k-1},\eta,x_{k+1},\ldots,x_n).$$

THEOREM 4.10. *Let $\alpha > \frac{1}{4}$ and $q > \frac{1}{\alpha}$. Then the second iterated integrals $\boldsymbol{\Gamma}(\varepsilon)_t - \boldsymbol{\Gamma}(\varepsilon)_s$, $\mathrm{Area}_{s,t}(\boldsymbol{\Gamma}(\varepsilon))$ converge with respect to the distance $d_1$, respectively, $d_2$.*

PROOF. The quantity one wants to estimate is the following:

(4.60)
$$\Delta(\varepsilon,\eta;t) = \mathbb{E}[(\mathrm{Area}_{0,t}(\Gamma(\varepsilon)) - \mathrm{Area}_{0,t}(\Gamma(\eta)))^2]$$
$$= \mathbb{E}\left[\left(\int_0^t \Gamma'(\varepsilon)_{x_1}^{(1)}\,dx_1 \int_0^{x_1} \Gamma'(\varepsilon)_{x_2}^{(2)}\,dx_2 \right.\right.$$
$$\left.\left. - \int_0^t \Gamma'(\eta)_{x_1}^{(1)}\,dx_1 \int_0^{x_1} \Gamma'(\eta)_{x_2}^{(2)}\,dx_2\right)^2\right]$$

for $\varepsilon > \eta > 0$, say.

By Corollary 4.8, it suffices to prove that, for $t \gg \varepsilon$, $\Delta(\varepsilon,\eta;t)$ is a finite sum of terms bounded by a constant times $t^{2\alpha_2}\varepsilon^{2\alpha_1}$ with $\alpha_1 > 0$ and $\alpha_1 + \alpha_2 = 2\alpha$. We shall nickname such a term an $(\alpha_1,\alpha_2)$-*term* for short. Note that an $(\alpha_1,\alpha_2)$-term is also necessarily an $(\alpha_1',\alpha_2')$-term for every $\alpha_1',\alpha_2'$ such that $\alpha_1' + \alpha_2' = 2\alpha$ and $0 < \alpha_1' < \alpha_1$ (our estimates are not necessarily optimal).

Now

$$\mathrm{Area}_{0,t}(\Gamma(\varepsilon)) - \mathrm{Area}_{0,t}(\Gamma(\eta))$$
$$= \int_0^t \Gamma'(\varepsilon)_{x_1}^{(1)}\,dx_1 \int_0^{x_1} (\Gamma'(\varepsilon)_{x_2}^{(2)} - \Gamma'(\eta)_{x_2}^{(2)})\,dx_2$$
$$+ \int_0^t (\Gamma'(\varepsilon)_{x_1}^{(1)} - \Gamma'(\eta)_{x_1}^{(1)})\,dx_1 \int_0^{x_1} \Gamma'(\eta)_{x_2}^{(2)}\,dx_2$$
(4.61)
$$= \int_0^t \Gamma'(\varepsilon)_{x_1}^{(1)}((\Gamma(\varepsilon)_{x_1}^{(2)} - \Gamma(\varepsilon)_0^{(2)}) - (\Gamma(\eta)_{x_1}^{(2)} - \Gamma(\eta)_0^{(2)}))\,dx_1$$
$$+ ((\Gamma(\varepsilon)_t^{(1)} - \Gamma(\varepsilon)_0^{(1)}) - (\Gamma(\eta)_t^{(1)} - \Gamma(\eta)_0^{(1)}))(\Gamma(\eta)_t^{(2)} - \Gamma(\eta)_0^{(2)})$$
$$- \int_0^t \Gamma'(\eta)_{x_1}^{(2)}\,dx_1((\Gamma(\varepsilon)_{x_1}^{(1)} - \Gamma(\varepsilon)_0^{(1)}) - (\Gamma(\eta)_{x_1}^{(1)} - \Gamma(\eta)_0^{(1)}))$$
$$= E_1(\varepsilon,\eta;t) + E_2(\varepsilon,\eta;t) + E_3(\varepsilon,\eta;t)$$

is the sum of three terms. Note that $E_3(\varepsilon,\eta;t)$ is of the same form as $E_1(\varepsilon,\eta;t)$, except for the permutation $\mathbf{\Gamma}^{(1)} \leftrightarrow \mathbf{\Gamma}^{(2)}$ and $\varepsilon \leftrightarrow \eta$. Let us estimate $\mathbb{E}[E_1(\varepsilon,\eta;t)^2]$ and $\mathbb{E}[E_2(\varepsilon,\eta;t)^2]$ separately.

*Estimation of $E_2$.* The first and second component of $\mathbf{\Gamma}$ decouple and one is left with

$$\mathbb{E}[E_2(\varepsilon,\eta;t)^2] = \mathbb{E}[((\Gamma(\varepsilon)_t^{(1)} - \Gamma(\eta)_t^{(1)}) - (\Gamma(\varepsilon)_0^{(1)} - \Gamma(\eta)_0^{(1)}))^2]$$
(4.62)
$$\times \mathbb{E}[(\Gamma(\eta)_t^{(2)} - \Gamma(\eta)_0^{(2)})^2]$$
$$\leq C|\varepsilon - \eta|^{2\alpha} t^{2\alpha}$$

by Lemma 1.5, which contributes an $(\alpha,\alpha)$-term to the estimate of Corollary 4.8, as explained at the very beginning of this proof.

*Estimation of $E_1$.* By rewriting $(\Gamma(\varepsilon)_{x_1}^{(2)} - \Gamma(\varepsilon)_0^{(2)}) - (\Gamma(\eta)_{x_1}^{(2)} - \Gamma(\eta)_0^{(2)})$ as

$$(\Gamma(\varepsilon)_{x_1}^{(2)} - \Gamma(\eta)_{x_1}^{(2)}) - (\Gamma(\varepsilon)_0^{(2)} - \Gamma(\eta)_0^{(2)})$$
$$= 2\operatorname{Re}\left(\int_{x_1+i\eta}^{x_1+i\varepsilon} \Gamma'^{+,(2)}(z)\,dz - \int_{i\eta}^{i\varepsilon} \Gamma'^{+,(2)}(z)\,dz\right),$$

one gets

(4.63) $$\mathbb{E}[E_1(\varepsilon,\eta;t)^2] = \mathcal{W}_1(\varepsilon,\eta)_t + \mathcal{W}_2(\varepsilon,\eta)_t,$$

where

$$\mathcal{W}_1(\varepsilon,\eta)_t = 2\operatorname{Re}\int_0^t dx_1 \int_0^{\varepsilon-\eta} d\xi_1 \int_0^t dy_1$$
$$\times \int_0^{\varepsilon-\eta} d\xi_2 (-i(x_1-y_1)+2\varepsilon)^{2\alpha-2}$$
$$\times ((-i(x_1-y_1)+\xi_1+\xi_2+2\eta)^{2\alpha-2}$$
$$\quad - (-ix_1+\xi_1+\xi_2+2\eta)^{2\alpha-2}$$
$$\quad - (iy_1+\xi_1+\xi_2+2\eta)^{2\alpha-2} + (\xi_1+\xi_2+2\eta)^{2\alpha-2})$$
$$= \frac{1}{2\alpha(2\alpha-1)} \cdot 2\operatorname{Re}\int_0^t dx_1 \int_0^t dx_2 (-i(x_1-y_1)+2\varepsilon)^{2\alpha-2}$$
$$\times [(-i(x_1-y_1)+2\eta)^{2\alpha} - 2(-i(x_1-y_1)+\varepsilon+\eta)^{2\alpha}$$
$$\quad + (-i(x_1-y_1)+2\varepsilon)^{2\alpha}]$$
(4.64)
$$\quad - \frac{1}{2\alpha(2\alpha-1)} \cdot 2\operatorname{Re}\int_0^t dx_1 \int_0^t dx_2 (-i(x_1-y_1)+2\varepsilon)^{2\alpha-2}$$

A STOCHASTIC CALCULUS BY ANALYTIC EXTENSION 39

$$\times [(-ix_1 + 2\eta)^{2\alpha} - 2(-ix_1 + \varepsilon + \eta)^{2\alpha} + (-ix_1 + 2\varepsilon)^{2\alpha}]$$

$$- \frac{1}{2\alpha(2\alpha - 1)} \cdot 2\operatorname{Re} \int_0^t dx_1 \int_0^t dx_2 (-i(x_1 - y_1) + 2\varepsilon)^{2\alpha - 2}$$

$$\times [(iy_1 + 2\eta)^{2\alpha} - 2(iy_1 + \varepsilon + \eta)^{2\alpha} + (iy_1 + 2\varepsilon)^{2\alpha}]$$

$$+ \frac{1}{2\alpha(2\alpha - 1)} \cdot 2\operatorname{Re} \int_0^t dx_1 \int_0^t dx_2 (-i(x_1 - y_1) + 2\varepsilon)^{2\alpha - 2}$$

$$\times [(2\eta)^{2\alpha} - 2(\varepsilon + \eta)^{2\alpha} + (2\varepsilon)^{2\alpha}]$$

and

$$\begin{aligned}
\mathcal{W}_2(\varepsilon, \eta)_t &= 2\operatorname{Re} \int_0^t dx_1 \int_0^{\varepsilon-\eta} d\xi_1 \int_0^t dy_1 \int_0^{\varepsilon-\eta} d\xi_2 (-i(x_1 - y_1) + 2\varepsilon)^{2\alpha - 2} \\
&\quad \times ((i(x_1 - y_1) + \xi_1 + \xi_2 + 2\eta)^{2\alpha - 2} \\
&\quad - (ix_1 + \xi_1 + \xi_2 + 2\eta)^{2\alpha - 2} \\
&\quad - (-iy_1 + \xi_1 + \xi_2 + 2\eta)^{2\alpha - 2} + (\xi_1 + \xi_2 + 2\eta)^{2\alpha - 2}) \\
&= \frac{1}{2\alpha(2\alpha - 1)} \cdot 2\operatorname{Re} \int_0^t dx_1 \int_0^t dx_2 (-i(x_1 - y_1) + 2\varepsilon)^{2\alpha - 2} \\
&\quad \times [(i(x_1 - y_1) + 2\eta)^{2\alpha} - 2(i(x_1 - y_1) + \varepsilon + \eta)^{2\alpha} \\
&\quad\quad\quad\quad\quad\quad + (i(x_1 - y_1) + 2\varepsilon)^{2\alpha}] \\
&\quad - \frac{1}{2\alpha(2\alpha - 1)} \cdot 2\operatorname{Re} \int_0^t dx_1 \int_0^t dx_2 (-i(x_1 - y_1) + 2\varepsilon)^{2\alpha - 2} \\
&\quad \times [(ix_1 + 2\eta)^{2\alpha} - 2(ix_1 + \varepsilon + \eta)^{2\alpha} + (ix_1 + 2\varepsilon)^{2\alpha}] \\
&\quad - \frac{1}{2\alpha(2\alpha - 1)} \cdot 2\operatorname{Re} \int_0^t dx_1 \int_0^t dx_2 (-i(x_1 - y_1) + 2\varepsilon)^{2\alpha - 2} \\
&\quad \times [(-iy_1 + 2\eta)^{2\alpha} - 2(-iy_1 + \varepsilon + \eta)^{2\alpha} + (-iy_1 + 2\varepsilon)^{2\alpha}] \\
&\quad + \frac{1}{2\alpha(2\alpha - 1)} \cdot 2\operatorname{Re} \int_0^t dx_1 \int_0^t dx_2 (-i(x_1 - y_1) + 2\varepsilon)^{2\alpha - 2} \\
&\quad \times [(2\eta)^{2\alpha} - 2(\varepsilon + \eta)^{2\alpha} + (2\varepsilon)^{2\alpha}].
\end{aligned} \qquad (4.65)$$

One has found

$$\mathcal{W}_1(\varepsilon, \eta)_t = \mathcal{V}_1(\varepsilon, \varepsilon)_t - 2\mathcal{V}_1\left(\varepsilon, \frac{\varepsilon + \eta}{2}\right)_t + \mathcal{V}_1(\varepsilon, \eta)_t = \Delta_2(\mathcal{V}_1(\varepsilon, \cdot))(\varepsilon, \eta)$$

(see Definition 4.9 above) and a similar relation for $\mathcal{W}_2(\varepsilon, \eta)_t$ in terms of the second variation of $\mathcal{V}_2$ (see proof of Theorem 4.4 for the definition of $\mathcal{V}_1$ and $\mathcal{V}_2$). Note that the second variation is at most of the same order as the first variation (i.e., the difference between the values at $\varepsilon$ and $\eta$)—for *regular* functions, it is even of strictly lower order [this is not necessarily true otherwise, take, e.g., $f(\varepsilon) = \varepsilon^\beta$, $\beta > 0$ small!]. Yet all bounds below are obtained by estimating the first variation, which is enough for our purposes.

Let us write asymptotic expansions for $\varepsilon, \eta \to 0$, $\varepsilon > \eta$, $\frac{\varepsilon}{\eta} \to 1$, $\varepsilon, \eta \ll t$ of each of the terms found in the computation of $\mathcal{V}_1$ and $\mathcal{V}_2$ in Theorem 4.4.

First,

$$\operatorname{Re} I_1(0, 0; 2\alpha - 2, 2\alpha; \varepsilon, \eta; 0, t)$$

$$(4.66) \quad = \frac{1}{4\alpha - 1} \operatorname{Re} i \left\{ (-it)^{4\alpha - 1} \left( 1 + 2i(4\alpha - 1)\frac{\eta}{t} + O\left(\frac{\eta^2}{t^2}\right) \right) \right.$$

$$\left. \times \left( 1 - \frac{(2 - 2\alpha)(1 - 4\alpha)}{2 - 4\alpha} \frac{2(\varepsilon - \eta)}{-it} + O\left(\frac{\eta^2}{t^2}\right) \right) \right\},$$

hence,

$$|\Delta_2(t \operatorname{Re} I_1(0, 0; 2\alpha - 2, 2\alpha; \varepsilon, \cdot; 0, t))(\varepsilon, \eta)| \leq C \eta t^{4\alpha - 1}$$

and similarly for

$$|\Delta_2(t \operatorname{Re} I_1(0, 0; 2\alpha - 1, 2\alpha; \varepsilon, \cdot; 0, t))(\varepsilon, \eta)|,$$

since the arguments of the hypergeometric functions (see 4.31) keep away from 1. Next,

$$|\Delta_2((\varepsilon + \eta) I_1(0, 0; 2\alpha - 1, 2\alpha - 1; \varepsilon, \cdot; 0, t))(\varepsilon, \eta)| \leq C(\eta t^{4\alpha - 1} + \eta^{4\alpha});$$

$$|\Delta_2(\eta \to (-it + 2\varepsilon)^{2\alpha - 1} t(-it + 2\eta)^{2\alpha})(\varepsilon, \eta)| \leq C \eta t^{4\alpha - 1}$$

and similarly for $|\Delta_2(\eta \to ((-it + 2\varepsilon)(-it + 2\eta))^{2\alpha})(\varepsilon, \eta)|$ and

$$\Delta_2(\eta \to (-it + 2\varepsilon)^{2\alpha - 1} t(-it + 2\eta)^{2\alpha})(\varepsilon, \eta);$$

$$|\Delta_2(\eta \to (\varepsilon \eta)^{2\alpha})(\varepsilon, \eta)| \leq C \eta^{4\alpha}.$$

In order to estimate $I_1(t, 0; 2\alpha - 1, 2\alpha; \varepsilon, \eta; 0, t)$, we apply the connection formula (4.4)

$$(4.67) \quad \begin{aligned} {}_2F_1\left(1 - 2\alpha, 1 + 2\alpha; 2 + 2\alpha; -\left(\frac{2\eta - it}{2(\varepsilon - \eta) + it}\right)\right) \\ = \frac{\Gamma(2 + 2\alpha)\Gamma(2\alpha)}{\Gamma(1 + 4\alpha)} {}_2F_1\left(1 - 2\alpha, 1 + 2\alpha; 1 - 2\alpha; \frac{2\varepsilon}{2(\varepsilon - \eta) + it}\right) \\ + \left(\frac{2\varepsilon}{2(\varepsilon - \eta) + it}\right)^{2\alpha} \frac{\Gamma(2 + 2\alpha)\Gamma(-2\alpha)}{\Gamma(1 + 2\alpha)\Gamma(1 - 2\alpha)} \\ \times {}_2F_1\left(1 + 4\alpha, 1; 1 + 2\alpha; \frac{2\varepsilon}{2(\varepsilon - \eta) + it}\right). \end{aligned}$$

and find out that the term on the first line of (4.34) has a second variation bounded by $C(t^{4\alpha-1}\eta + t^{2\alpha}\eta^{2\alpha})$, while the second variation of the term on the second line is bounded by $Ct^{2\alpha-1}\eta^{2\alpha+1}$.

Finally, $\mathcal{V}_{1,4}$ is an $(\alpha, \alpha)$-term.

On the whole, we have proved that $\mathcal{W}_1(\varepsilon, \eta)_t$ is a sum of $(\alpha_1, \alpha_2)$-terms for $\alpha_1 = \frac{1}{2}, \alpha, 2\alpha$ or $\alpha + \frac{1}{2}$.

The same goes for all terms in $\mathcal{W}_2(\varepsilon, \eta)_t$, except of course for the singular term in $tI_2(0, 0; 2\alpha - 2, 2\alpha; \varepsilon, \eta; 0, t)$, namely,

$$f(\varepsilon, \eta) = -\frac{1}{4\alpha - 1} i e^{2i\pi\alpha} t(2\eta)^{4\alpha-1}$$

$$\times \lim_{s \overset{>}{\to} 0} {}_2F_1\left(2 - 2\alpha, 1 - 4\alpha; 2 - 4\alpha; \frac{2(\varepsilon + \eta)}{2\eta - is}\right);$$

see the second line of (4.38).

The connection formula (4.3) yields

(4.68)
$$\lim_{s \overset{>}{\to} 0} {}_2F_1\left(2 - 2\alpha, 1 - 4\alpha; 2 - 4\alpha; \frac{2(\varepsilon + \eta)}{2\eta - is}\right)$$
$$= \Gamma(2 - 4\alpha)$$
$$\times \left\{ \frac{\Gamma(-1 - 2\alpha)}{\Gamma(1 - 4\alpha)\Gamma(-2\alpha)} \left(\frac{\varepsilon + \eta}{\eta}\right)^{2\alpha - 2} e^{-2i\pi\alpha} \right.$$
$$\times {}_2F_1\left(2 - 2\alpha, 1 + 2\alpha; 2 + 2\alpha; \frac{\eta}{\varepsilon + \eta}\right)$$
$$\left. + \frac{\Gamma(1 + 2\alpha)}{\Gamma(2 - 2\alpha)} \left(\frac{\varepsilon + \eta}{\eta}\right)^{4\alpha - 1} e^{-i\pi(4\alpha - 1)} \right\}.$$

Suppose $\alpha > \frac{1}{4}$: then quite simply

$$|\Delta_2 f(\varepsilon, \eta)| \leq C\eta^{4\alpha-1} t,$$

hence, $f(\varepsilon, \eta)$ is a $(\frac{4\alpha-1}{2}, \frac{1}{2})$-term. $\square$

4.3. *Convergence in the $d_3$-distance for $\frac{1}{4} < \alpha \leq \frac{1}{3}$.* An analogue of Proposition 4.5 holds for the third chaos of FBM, namely:

PROPOSITION 4.11 (See [5]). *Let $\mathbf{w}, \mathbf{v}$ be two $V$-valued paths with $q$-bounded variation and bounded $d_3$-norm. Then for every $\kappa > q/3 - 1$, there exists a constant $C$ depending only on $q$ and $\kappa$ such that*

$$(d_3(\mathbf{w}, \mathbf{v})_{s,t})^{q/3} \leq C \sum_{n \geq 1} n^\kappa \sum_{k=1}^{2^n} |\mathbf{w}^3_{t^n_{k-1}, t^n_k} - \mathbf{v}^3_{t^n_{k-1}, t^n_k}|^{q/3}$$

$$+ C \left( \sum_{n \geq 1} n^{\kappa} \sum_{k=1}^{2^n} |\mathbf{w}^2_{t^n_{k-1},t^n_k} - \mathbf{v}^2_{t^n_{k-1},t^n_k}|^{q/2} \right)^{2/3}$$

$$\times \left( \sum_{n \geq 1} n^{\kappa} \sum_{k=1}^{2^n} |\mathbf{w}^1_{t^n_{k-1},t^n_k}|^q + |\mathbf{v}^1_{t^n_{k-1},t^n_k}|^q \right)^{1/3}$$

(4.69)
$$+ C \left( \sum_{n \geq 1} n^{\kappa} \sum_{k=1}^{2^n} |\mathbf{w}^1_{t^n_{k-1},t^n_k} - \mathbf{v}^1_{t^n_{k-1},t^n_k}|^q \right)^{1/3}$$

$$\times \left( \sum_{n \geq 1} n^{\kappa} \sum_{k=1}^{2^n} |\mathbf{w}^2_{t^n_{k-1},t^n_k}|^{q/2} + |\mathbf{v}^2_{t^n_{k-1},t^n_k}|^{q/2} \right)^{2/3}$$

$$+ C \left( \sum_{n \geq 1} n^{\kappa} \sum_{k=1}^{2^n} |\mathbf{w}^1_{t^n_{k-1},t^n_k} - \mathbf{v}^1_{t^n_{k-1},t^n_k}|^q \right)^{1/3}$$

$$\times \left( \sum_{n \geq 1} n^{\kappa} \sum_{k=1}^{2^n} |\mathbf{w}^1_{t^n_{k-1},t^n_k}|^q + |\mathbf{v}^1_{t^n_{k-1},t^n_k}|^q \right)^{2/3}.$$

DEFINITION 4.12. Let

(4.70)
$$\mathrm{Vol}_{0,t}(\varepsilon_1, \varepsilon_2, \varepsilon_3)$$
$$:= \int_0^t \Gamma'(\varepsilon_1)^{(1)}_{x_1} \, dx_1 \int_0^{x_1} \Gamma'(\varepsilon_2)^{(2)}_{x_2} \, dx_2$$
$$\times \int_0^{x_2} \Gamma'(\varepsilon_3)^{(3)}_{x_3} \, dx_3.$$

The volume functional generated by the $\varepsilon$-approximation of a three-dimensional fractional Brownian motion is given by $\mathrm{Vol}(\varepsilon, \varepsilon, \varepsilon)$, which is more appropriately written as $\mathrm{Vol}(\mathbf{\Gamma}(\varepsilon))$.

THEOREM 4.13. *Let $\alpha > \frac{1}{4}$ and $q > \frac{1}{\alpha}$. Then the third iterated integrals $\mathrm{Vol}_{s,t}(\mathbf{\Gamma}(\varepsilon))$ converge with respect to the distance $d_3$.*

The rest of the article is dedicated to the proof of this theorem.

DEFINITION 4.14. Let

(4.71)
$$\mathcal{W}(\varepsilon_1, \varepsilon_2, \varepsilon_3)_t := \mathbb{E}[\mathrm{Vol}_{0,t}(\varepsilon_1, \varepsilon_2, \varepsilon_3)^2]$$
$$= \sum_{\boldsymbol{\sigma}} \int_0^t dx_1 \int_0^t dy_1 \cdot \int_0^{x_1} dx_2 \int_0^{y_1} dy_2 \cdot \int_0^{x_2} dx_3$$

$$\times \int_0^{y_2} dy_3 \prod_j [-i\sigma_j(x_j - y_j) + 2\varepsilon_j]^{2\alpha-2},$$

where $\boldsymbol{\sigma} = (\sigma_1, \sigma_2, \sigma_3) \in \{\pm 1\}^3$.

Rewrite $\mathrm{Vol}_{0,t}(\boldsymbol{\Gamma}(\varepsilon)) - \mathrm{Vol}_{0,t}(\boldsymbol{\Gamma}(\eta))$ as

$$[\mathrm{Vol}_{0,t}(\varepsilon, \varepsilon, \varepsilon) - \mathrm{Vol}_{0,t}(\varepsilon, \varepsilon, \eta)] + [\mathrm{Vol}_{0,t}(\varepsilon, \varepsilon, \eta) - \mathrm{Vol}_{0,t}(\varepsilon, \eta, \eta)]$$
$$+ [\mathrm{Vol}_{0,t}(\varepsilon, \eta, \eta) - \mathrm{Vol}_{0,t}(\eta, \eta, \eta)].$$

Then

(4.72)
$$\mathbb{E}[(\mathrm{Vol}_{0,t}(\boldsymbol{\Gamma}(\varepsilon)) - \mathrm{Vol}_{0,t}(\boldsymbol{\Gamma}(\eta)))^2]$$
$$\leq C(\Delta_2(\mathcal{W}(\varepsilon, \varepsilon, \cdot))(\varepsilon, \eta) + \Delta_2(\mathcal{W}(\varepsilon, \cdot, \eta))(\varepsilon, \eta)$$
$$+ \Delta_2(\mathcal{W}(\cdot, \eta, \eta))(\varepsilon, \eta)).$$

Corollary 4.8 may be restated in this case: convergence in the $d_3$-norm holds true if the following estimate holds:

(4.73)
$$\mathbb{E}[\|\boldsymbol{\Gamma}(\varepsilon)_{s,t}^2 - \boldsymbol{\Gamma}(\eta)_{s,t}^2\|^2] \leq \sum_{i=1}^I C_i \varepsilon^{3\alpha_{i,1}} |t-s|^{3\alpha_{i,2}}$$

for some constants $C_i > 0$ and couples of exponents $(\alpha_{i_1}, \alpha_{i_2})$ such that $\alpha_{i_1} + \alpha_{i_2} = 2\alpha$ and $\alpha_{i_1} > 0$.

In the present case (cf. with the definition given in the course of the proof of Theorem 4.10), a function $f(t, \varepsilon, \eta)$ will be called an $(\alpha_1, \alpha_2)$-*term* $(\alpha_1 > 0, \alpha_1 + \alpha_2 = 2\alpha)$ if

$$|f(t, \varepsilon, \eta)| \leq C\varepsilon^{3\alpha_1} t^{3\alpha_2}.$$

So one is left with the problem of giving $(\alpha_1, \alpha_2)$-type estimates for the above $\Delta_2$-terms with $\varepsilon = j^{-\beta}$, $\eta = (j+1)^{-\beta}$ for some $\beta$ (so that $\varepsilon, \eta \to 0$ and $\frac{\varepsilon}{\eta} \to 1$). In the following estimates of $\mathcal{W}(\varepsilon_1, \varepsilon_2, \varepsilon_3)_t$, one may (and shall!) assume that $\varepsilon_1, \varepsilon_2, \varepsilon_3$ are of the same order, $\frac{\varepsilon_1}{\varepsilon_2} \simeq \frac{\varepsilon_1}{\varepsilon_3} \simeq 1$, and that $t \gg \varepsilon_1, \varepsilon_2, \varepsilon_3$. We shall sometimes write $\varepsilon = \sup(\varepsilon_1, \varepsilon_2, \varepsilon_3)$ for short.

The explicit integration

(4.74)
$$\int_0^{x_2} dx_3 \int_0^{y_2} dy_3 (-i\sigma_3(x_3 - y_3) + 2\varepsilon_3)^{2\alpha-2}$$
$$= ((2\varepsilon_3)^{2\alpha} - (-i\sigma_3 x_2 + 2\varepsilon_3)^{2\alpha} - (i\sigma_3 y_2 + 2\varepsilon_3)^{2\alpha}$$
$$+ (-i\sigma_3(x_2 - y_2) + 2\varepsilon_3)^{2\alpha})(2\alpha(2\alpha - 1))^{-1}$$

yields $\mathcal{W}(\varepsilon_1, \varepsilon_2, \varepsilon_3)_t = \sum_{j=1}^4 \mathcal{W}_j(\varepsilon_1, \varepsilon_2, \varepsilon_3)_t$. We shall estimate these four terms separately.

*Estimation of* $\mathcal{W}_1$. Since

$$\mathcal{W}_1(\varepsilon_1,\varepsilon_2,\varepsilon_3)_t = \frac{(2\varepsilon_3)^{2\alpha}}{2\alpha(2\alpha-1)}\mathcal{V}(\varepsilon_1,\varepsilon_2)_t,$$

this term goes to 0 when $\varepsilon_1,\varepsilon_2,\varepsilon_3 \to 0$ if $\alpha > \frac{1}{4}$ and one easily gets $(\alpha_1,\alpha_2)$-type estimates by re-employing those obtained beforehand for Lévy's area.

*Estimation of* $\mathcal{W}_2$. Two successive integrations by parts yield

$$\mathcal{W}_2(\varepsilon_1,\varepsilon_2,\varepsilon_3)_t$$
$$= \sum_\sigma \frac{i\sigma_2}{2\alpha(2\alpha-1)^2} \int_0^t dx_1 \int_0^t dy_1 (-i\sigma_1(x_1-y_1)+2\varepsilon_1)^{2\alpha-2}$$
$$\times \int_0^{x_1} dx_2(-i\sigma_3 x_2+2\varepsilon_3)^{2\alpha}$$
$$\times [(-i\sigma_2(x_2-y_1)+2\varepsilon_2)^{2\alpha-1} - (-i\sigma_2 x_2+2\varepsilon_2)^{2\alpha-1}]$$
$$= -\sum_\sigma \frac{\sigma_1\sigma_2}{2\alpha(2\alpha-1)^3}$$
$$\times \Bigg[ \int_0^t dy_1(-i\sigma_1(t-y_1)+2\varepsilon_1)^{2\alpha-1}$$

(4.75)
$$\times \int_0^t dx_2(-i\sigma_3 x_2+2\varepsilon_3)^{2\alpha}(-i\sigma_2(x_2-y_1)+2\varepsilon_2)^{2\alpha-1}$$
$$- \int_0^t dx_1(-i\sigma_3 x_1+2\varepsilon_3)^{2\alpha}$$
$$\times \int_0^t dy_1(-i\sigma_1(x_1-y_1)+2\varepsilon_1)^{2\alpha-1}$$
$$\times (-i\sigma_2(x_1-y_1)+2\varepsilon_2)^{2\alpha-1}$$
$$- \int_0^t dx_1[(-i\sigma_1(x_1-t)+2\varepsilon_1)^{2\alpha-1}$$
$$- (-i\sigma_1 x_1+2\varepsilon_1)^{2\alpha-1}]$$
$$\times \int_0^{x_1} dx_2(-i\sigma_3 x_2+2\varepsilon_3)^{2\alpha}(-i\sigma_2 x_2+2\varepsilon_2)^{2\alpha-1}\Bigg],$$

which is the sum of three terms, $\mathcal{W}_{2,i}$, $i=1,2,3$.

(i) (Estimation of $\mathcal{W}_{2,1}$.)
By Lemma 4.1,

$$\int_0^t dx_2(-i\sigma_3 x_2+2\varepsilon_3)^{2\alpha}(-i\sigma_2(x_2-y_1)+2\varepsilon_2)^{2\alpha-1}$$

$$
\begin{aligned}
&= \delta_{\sigma_3,1}\delta_{\sigma_2,1} I_1(0,y_1;2\alpha,2\alpha-1;\varepsilon_3,\varepsilon_2;0,t)\\
&\quad + \delta_{\sigma_3,-1}\delta_{\sigma_2,-1}\overline{I_1(0,y_1;2\alpha,2\alpha-1;\varepsilon_3,\varepsilon_2;0,t)}\\
&\quad + \delta_{\sigma_3,-1}\delta_{\sigma_2,1} I_2(0,y_1;2\alpha,2\alpha-1;\varepsilon_3,\varepsilon_2;0,t)\\
&\quad + \delta_{\sigma_3,1}\delta_{\sigma_2,-1}\overline{I_2(0,y_1;2\alpha,2\alpha-1;\varepsilon_3,\varepsilon_2;0,t)}.
\end{aligned}
\tag{4.76}
$$

If $\sigma_2 = \sigma_3 = 1$ and $\varepsilon_3 > \varepsilon_2$, then $I_1(0,y_1;2\alpha,2\alpha-1;\varepsilon_3,\varepsilon_2;0,t) = F_1(0,y_1;2\alpha,2\alpha-1;\varepsilon_3,\varepsilon_2;t) - F_1(0,y_1;2\alpha,2\alpha-1;\varepsilon_3,\varepsilon_2;0)$ (see Lemma 4.1), with

$$
\begin{aligned}
&F_1(0,y_1;2\alpha,2\alpha-1;\varepsilon_3,\varepsilon_2;t)\\
&= i\frac{(2\varepsilon_2 - i(t-y_1))^{2\alpha}}{2\alpha}(2(\varepsilon_3-\varepsilon_2) - iy_1)^{2\alpha}\\
&\quad \times {}_2F_1\left(-2\alpha, 2\alpha; 2\alpha+1; -\frac{2\varepsilon_2 - i(t-y_1)}{2(\varepsilon_3-\varepsilon_2) - iy_1}\right).
\end{aligned}
\tag{4.77}
$$

The real part of the argument $z(t) = -\frac{2\varepsilon_2 - i(t-y_1)}{2(\varepsilon_3-\varepsilon_2) - iy_1}$ of the hypergeometric function is negative, hence, $F_1(0,y_1;2\alpha,2\alpha-1;\varepsilon_3,\varepsilon_2;t)$ is regular in $\varepsilon_3,\varepsilon_2,t$. In the region $\frac{y_1}{t} \simeq 1$, one has $|z(t)| \ll 1$, hence (by Taylor's formula), the associated second variations are less than a constant times

$$
\begin{aligned}
&\int_0^t dy_1 ((t-y_1) + \varepsilon)^{2\alpha-1}\\
&\quad \times (((t-y_1)+\varepsilon)^{2\alpha-1}(y_1+\varepsilon)^{2\alpha} + ((t-y_1)+\varepsilon)^{2\alpha}(y_1+\varepsilon)^{2\alpha-1}),
\end{aligned}
\tag{4.78}
$$

which has a well defined limit for $\varepsilon \to 0$ as soon as $\alpha > \frac{1}{6}$, and one gets a $(\frac{1}{3}, 2\alpha - \frac{1}{3})$-bound for this term.

If, on the contrary, $\frac{y_1}{t} \ll 1$, then $|z(t)|$ is large and one must apply the first connection formula (4.3),

$$
\begin{aligned}
&{}_2F_1(-2\alpha, 2\alpha; 2\alpha+1; z(t))\\
&= \Gamma(2\alpha+1)\Bigg\{\frac{\Gamma(4\alpha)}{\Gamma(2\alpha)\Gamma(1+4\alpha)}\left(\frac{2\varepsilon_2 - i(t-y_1)}{2(\varepsilon_3-\varepsilon_2) - iy_1}\right)^{2\alpha}\\
&\quad \times {}_2F_1\left(-2\alpha, -4\alpha; 1-4\alpha; \frac{1}{z(t)}\right)\\
&\quad + \frac{\Gamma(-4\alpha)}{\Gamma(-2\alpha)}\left(\frac{2(\varepsilon_3-\varepsilon_2) - iy_1}{2\varepsilon_2 - i(t-y_1)}\right)^{2\alpha}\Bigg\}.
\end{aligned}
\tag{4.79}
$$

The large factor $(2(\varepsilon_3-\varepsilon_2) - iy_1)^{2\alpha}$ in the first term of the r.h.s. is compensated by the overall factor in $F_1$, hence, the same method as before gives a $(\frac{1}{3}, 2\alpha - \frac{1}{3})$-bound for this term $(\alpha > \frac{1}{6})$. The same goes for the second term of the r.h.s.

Now

$$
\begin{aligned}
F_1(0, y_1; 2\alpha, 2\alpha - 1; \varepsilon_3, \varepsilon_2; 0) \\
= i\frac{(2\varepsilon_2 + iy_1)^{2\alpha}}{2\alpha}(2(\varepsilon_3 - \varepsilon_2) - iy_1)^{2\alpha} \\
\times {}_2F_1\left(-2\alpha, 2\alpha; 2\alpha + 1; -\frac{2\varepsilon_2 + iy_1}{2(\varepsilon_3 - \varepsilon_2) - iy_1}\right)
\end{aligned}
\tag{4.80}
$$

can be estimated in the same way in the region $y_1 = O(\varepsilon)$ [note that $\operatorname{Re} z(0)$ is not necessarily positive, but $|\operatorname{Im} z(0)|$ is bounded from below]. In the region $y_1 \gg \varepsilon$ where $z(0) \simeq 1$, the hypergeometric function must be transformed by means of the second connection formula (4.4),

$$
\begin{aligned}
{}_2F_1\left(-2\alpha, 2\alpha; 2\alpha + 1; -\frac{2\varepsilon_2 + iy_1}{2(\varepsilon_3 - \varepsilon_2) - iy_1}\right) \\
= \Gamma(2\alpha + 1)\left\{\frac{\Gamma(2\alpha + 1)}{\Gamma(4\alpha + 1)} {}_2F_1\left(-2\alpha, 2\alpha; -2\alpha; \frac{2\varepsilon_3}{2(\varepsilon_3 - \varepsilon_2) - iy_1}\right)\right. \\
+ \frac{2\varepsilon_3}{2(\varepsilon_3 - \varepsilon_2) - iy_1} \frac{\Gamma(-2\alpha - 1)}{\Gamma(-2\alpha)\Gamma(2\alpha)} \\
\left.\times {}_2F_1\left(4\alpha + 1, 1; 2\alpha + 2; \frac{2\varepsilon_3}{2(\varepsilon_3 - \varepsilon_2) - iy_1}\right)\right\},
\end{aligned}
\tag{4.81}
$$

which can be estimated by means of Taylor's formula as before, with the same result.

Similar arguments hold if $\sigma_3 = -1, \sigma_2 = 1$; then $I_2(0, y_1; 2\alpha, 2\alpha-1; \varepsilon_3, \varepsilon_2; 0, t) = F_2(0, y_1; 2\alpha, 2\alpha - 1; \varepsilon_3, \varepsilon_2; t) - F_2(0, y_1; 2\alpha, 2\alpha - 1; \varepsilon_3, \varepsilon_2; 0)$ (see Lemma 4.1) with $F_2$ involving the same hypergeometric function but with a different argument, namely, $w(t) = \frac{2\varepsilon_2 - i(t-y_1)}{2(\varepsilon_3 + \varepsilon_2) + iy_1}$. The real part of $w(t)$ does not have a constant sign this time but keeps away from the forbidden line $[1; +\infty[$ (away from the real axis, actually). Details are left to the reader.

Note that the case $\varepsilon_3 < \varepsilon_2$ can be settled by the same type of arguments, for instance, by using an integration by parts to relate to the previous case.

(ii) (Estimation of $\mathcal{W}_{2,2}$.)

By Lemma 4.1,

$$
\begin{aligned}
\int_0^t dy_1(-i\sigma_1(x_1 - y_1) + 2\varepsilon_1)^{2\alpha-1}(-i\sigma_2(x_1 - y_1) + 2\varepsilon_2)^{2\alpha-1} \\
= \delta_{\sigma_1,-1}\delta_{\sigma_2,-1}I_1(0, 0; 2\alpha - 1, 2\alpha - 1; \varepsilon_1, \varepsilon_2; -x_1, t - x_1) \\
+ \delta_{\sigma_1,1}\delta_{\sigma_2,1}\overline{I_1(0, 0; 2\alpha - 1, 2\alpha - 1; \varepsilon_1, \varepsilon_2; -x_1, t - x_1)}
\end{aligned}
\tag{4.82}
$$

$$+ \delta_{\sigma_1,1}\delta_{\sigma_2,-1} I_2(0,0;2\alpha-1,2\alpha-1;\varepsilon_1,\varepsilon_2;-x_1,t-x_1)$$
$$+ \delta_{\sigma_1,-1}\delta_{\sigma_2,1} \overline{I_2(0,0;2\alpha-1,2\alpha-1;\varepsilon_1,\varepsilon_2;-x_1,t-x_1)}.$$

In the case $\sigma_1 = \sigma_2 = -1$, one writes (see Lemma 4.1)

(4.83)
$$\begin{aligned} I_1(0,0;2\alpha-1,2\alpha-1;\varepsilon_1,\varepsilon_2;-x_1,t-x_1) \\ = \Phi_1(0,0;2\alpha-1,2\alpha-1;\varepsilon_1,\varepsilon_2;t-x_1) \\ - \Phi_1(0,0;2\alpha-1,2\alpha-1;\varepsilon_1,\varepsilon_2;-x_1), \end{aligned}$$

where

(4.84)
$$\Phi_1(0,0;2\alpha-1,2\alpha-1;\varepsilon_1,\varepsilon_2;z)$$
$$= i\frac{(2\varepsilon_2 - iz)^{4\alpha-1}}{4\alpha-1} {}_2F_1\left(1-2\alpha, 1-4\alpha; 2-4\alpha; -\frac{2(\varepsilon_1-\varepsilon_2)}{2\varepsilon_2 - iz}\right).$$

Since the argument of the hypergeometric function keeps small, this case is easily settled by using Taylor's formula at order one as before: the associated second variations are less than a constant times

(4.85)
$$\int_0^t dx_1 (x_1 + \varepsilon_1)^{2\alpha}[(x_1 + \varepsilon_2)^{4\alpha-2} + ((t-x_1) + \varepsilon_2)^{4\alpha-2}]$$
$$+ \int_0^t dx_1 (x_1 + \varepsilon_1)^{2\alpha-1}[(x_1 + \varepsilon_2)^{4\alpha-1} + ((t-x_1) + \varepsilon_2)^{4\alpha-1}].$$

Hence, one gets once more a $(\frac{1}{3}, 2\alpha - \frac{1}{3})$-bound for this term.

Suppose now $\sigma_1 = 1$, $\sigma_2 = -1$: then

(4.86)
$$\begin{aligned} I_2(0,0;2\alpha-1,2\alpha-1;\varepsilon_1,\varepsilon_2;-x_1,t-x_1) \\ = F_2(0,0;2\alpha-1,2\alpha-1;\varepsilon_1,\varepsilon_2;t-x_1) \\ - F_2(0,0;2\alpha-1,2\alpha-1;\varepsilon_1,\varepsilon_2;-x_1), \end{aligned}$$

where

(4.87)
$$F_2(0,0;2\alpha-1,2\alpha-1;\varepsilon_1,\varepsilon_2;z)$$
$$= i\frac{(2\varepsilon_2 - iz)^{2\alpha}}{2\alpha-1}(2(\varepsilon_1+\varepsilon_2))^{2\alpha-1}$$
$$\times {}_2F_1\left(1-2\alpha, 2\alpha; 2\alpha+1; \frac{2\varepsilon_2 - iz}{2(\varepsilon_1+\varepsilon_2)}\right).$$

If $z = O(\varepsilon)$, $z = x_1$ or $t - x_1$, then the argument of the hypergeometric function is bounded and stays away from $[1; +\infty[$, and

(4.88)
$$|(2\varepsilon_2 - iz)^{2\alpha}(2(\varepsilon_1+\varepsilon_2))^{2\alpha-1}| \leq C\varepsilon^{4\alpha-1},$$

hence, one gets (for $\alpha > \frac{1}{4}$ this time) a $(\frac{4\alpha-1}{3}, \frac{2\alpha+1}{3})$-bound. Otherwise the first connection formula (4.3) yields

$$
\begin{aligned}
{}_2F_1&\left(1-2\alpha, 2\alpha; 2\alpha+1; \frac{2\varepsilon_2 - iz}{2(\varepsilon_1+\varepsilon_2)}\right) \\
&= \Gamma(1+2\alpha)\Bigg\{ \frac{\Gamma(4\alpha-1)}{\Gamma(2\alpha)\Gamma(4\alpha)} \left(-\frac{2(\varepsilon_1+\varepsilon_2)}{2\varepsilon_2-iz}\right)^{1-2\alpha} \\
&\qquad\qquad \times {}_2F_1\left(1-2\alpha, 1-4\alpha; 2-4\alpha; \frac{2(\varepsilon_1+\varepsilon_2)}{2\varepsilon_2-iz}\right) \\
&\qquad\qquad + \frac{\Gamma(1-4\alpha)}{\Gamma(1-2\alpha)} \left(-\frac{2(\varepsilon_1+\varepsilon_2)}{2\varepsilon_2-iz}\right)^{2\alpha} \Bigg\}.
\end{aligned}
\tag{4.89}
$$

The first term of the r.h.s. is similar to that obtained in the case $\sigma_1 = \sigma_2 = -1$, while the second term is a $(\frac{4\alpha-1}{3}, \frac{2\alpha+1}{3})$-term if $\alpha > \frac{1}{4}$.

(iii) (Estimation of $\mathcal{W}_{2,3}$.)

According to the relative signs of $\sigma_2, \sigma_3$,

$$\int_0^{x_1} dx_2 (-i\sigma_3 x_2 + 2\varepsilon_3)^{2\alpha} (-i\sigma_2 x_2 + 2\varepsilon_2)^{2\alpha-1}$$

may be written (up to conjugacy) as

$$\Phi_1(0,0; 2\alpha, 2\alpha-1; \varepsilon_3, \varepsilon_2; x_1) - \Phi_1(0,0; 2\alpha, 2\alpha-1; \varepsilon_3, \varepsilon_2; 0)$$

or

$$\Phi_2(0,0; 2\alpha, 2\alpha-1; \varepsilon_3, \varepsilon_2; x_1) - \lim_{u\to 0, u>0} \Phi_2(0,0; 2\alpha, 2\alpha-1; \varepsilon_3, \varepsilon_2; u).$$

The argument of the hypergeometric function appearing in $\Phi_1$ remains bounded, with a real part $\ll 1$. The second case (involving $\Phi_2$) is slightly more complicated; the argument of the hypergeometric function, $z(s) = \frac{2(\varepsilon_1+\varepsilon_2)}{2\varepsilon_2-is}$, must be evaluated for $s = x_1$ and $s = 0$. If $s \geq C\varepsilon_2$, then $|z(s)|$ is bounded away from 1 and the hypergeometric function is regular and bounded. Otherwise one may use the second connnection formula (4.4) as in (4.68). Details are left to the reader.

*Estimation of $\mathcal{W}_3$.* It can be obtained from the previous case by simply permuting the $x$-coordinates with the $y$-coordinates.

*Estimation of $\mathcal{W}_4$.* One computes

$$
\begin{aligned}
\mathcal{W}_4&(\varepsilon_1, \varepsilon_2, \varepsilon_3)_t \\
&= \int_0^t dx_1 \int_0^t dy_1 (-i\sigma_1(x_1-y_1) + 2\varepsilon_1)^{2\alpha-2} \\
&\quad \times \int_0^{x_1} \int_0^{y_1} dy_2 (-i\sigma_2(x_2-y_2) + 2\varepsilon_2)^{2\alpha-2} (-i\sigma_3(x_2-y_2) + 2\varepsilon_3)^{2\alpha}
\end{aligned}
$$

$$= -\frac{i\sigma_1}{2\alpha - 1} \int_0^t dx_1 (-i\sigma_1(x_1 - t) + 2\varepsilon_1)^{2\alpha-1}$$

$$\times \int_0^{x_1} dx_2 \int_0^t dy_2 (-i\sigma_2(x_2 - y_2) + 2\varepsilon_2)^{2\alpha-2}$$

$$\times (-i\sigma_3(x_2 - y_2) + 2\varepsilon_3)^{2\alpha}$$

$$+ \frac{i\sigma_1}{2\alpha - 1} \int_0^t dx_1 \int_0^t dy_1 (-i\sigma_1(x_1 - y_1) + 2\varepsilon_1)^{2\alpha-1}$$

(4.90) $$\times \int_0^{x_1} dx_2 (-i\sigma_2(x_2 - y_1) + 2\varepsilon_2)^{2\alpha-2}(-i\sigma_3(x_2 - y_1) + 2\varepsilon_3)^{2\alpha}$$

$$= -\frac{1}{2\alpha(2\alpha - 1)} \int_0^t dx_1 (-i\sigma_1(x_1 - t) + 2\varepsilon_1)^{2\alpha}$$

$$\times \int_0^t dy_2 (-i\sigma_2(x_1 - y_2) + 2\varepsilon_2)^{2\alpha-2}(-i\sigma_3(x_1 - y_2) + 2\varepsilon_3)^{2\alpha}$$

$$- \frac{1}{2\alpha(2\alpha - 1)} \int_0^t dy_1 (-i\sigma_1(t - y_1) + 2\varepsilon_1)^{2\alpha}$$

$$\times \int_0^t dx_2 (-i\sigma_2(x_2 - y_1) + 2\varepsilon_2)^{2\alpha-2}(-i\sigma_3(x_2 - y_1) + 2\varepsilon_3)^{2\alpha}$$

$$+ \frac{1}{2\alpha(2\alpha - 1)} \int_0^t dx_1 \int_0^t dy_1 (-i\sigma_1(x_1 - y_1) + 2\varepsilon_1)^{2\alpha}$$

$$\times (-i\sigma_2(x_1 - y_1) + 2\varepsilon_2)^{2\alpha-2}$$

$$\times (-i\sigma_3(x_1 - y_1) + 2\varepsilon_3)^{2\alpha}.$$

The last integral may be rewritten as

$$\int_{-t}^t dx (-i\sigma_1 x + 2\varepsilon_1)^{2\alpha} (-i\sigma_2 x + 2\varepsilon_2)^{2\alpha-2} (-i\sigma_3 x + 2\varepsilon_3)^{2\alpha} (t - |x|).$$

Since

$$([(-i\sigma_1 x + 2\varepsilon_1)(-i\sigma_3 x + 2\varepsilon_3)]^{2\alpha+1})'$$
$$= (2\alpha + 1)[(-i\sigma_1 x + 2\varepsilon_1)(-i\sigma_3 x + 2\varepsilon_3)]^{2\alpha}[-2\sigma_1\sigma_3 x - 2i(\sigma_1\varepsilon_3 + \sigma_3\varepsilon_1)],$$

one gets

$$\int_0^t dx (-i\sigma_1 x + 2\varepsilon_1)^{2\alpha}(-i\sigma_2 x + 2\varepsilon_2)^{2\alpha-2}(-i\sigma_3 x + 2\varepsilon_3)^{2\alpha} x$$

$$= -\frac{i}{\sigma_1\sigma_3}(\sigma_1\varepsilon_3 + \sigma_3\varepsilon_1)$$

(4.91) $$\times \int_0^t dx (-i\sigma_1 x + 2\varepsilon_1)^{2\alpha}(-i\sigma_2 x + 2\varepsilon_2)^{2\alpha-2}(-i\sigma_3 x + 2\varepsilon_3)^{2\alpha}$$

$$-\frac{1}{2(2\alpha+1)\sigma_1\sigma_3}$$
$$\times \int_0^t dx(-i\sigma_2 x + 2\varepsilon_2)^{2\alpha-2}$$
$$\times [((-i\sigma_1 x + 2\varepsilon_1)(-i\sigma_3 x + 2\varepsilon_3))^{2\alpha+1}]'.$$

Finally, one is left with the problem of estimating the second variation of

$$J_{\beta_1,\beta_2,\beta_3} := \int_0^t dx(-i\sigma_1 x + 2\varepsilon_1)^{\beta_1}(-i\sigma_2 x + 2\varepsilon_2)^{\beta_2}(-i\sigma_3 x + 2\varepsilon_3)^{\beta_3}$$

with $(\beta_1,\beta_2,\beta_3) = (2\alpha, 2\alpha-2, 2\alpha)$ or $(2\alpha+1, 2\alpha-3, 2\alpha+1)$. Taylor's formula at order 1 yields very easily a $(\frac{1}{3}, 2\alpha - \frac{1}{3})$-bound for $\alpha > \frac{1}{6}$.

INSTITUT ELIE CARTAN
LABORATOIRE ASSOCIÉ AU CNRS UMR 7502
UNIVERSITÉ HENRI POINCARÉ NANCY I
B.P. 239
F-54506 VANDŒUVRE LÈS NANCY CEDEX
FRANCE
E-MAIL: jeremie.unterberger@iecn.u-nancy.fr